\documentclass[reqnox,10pt]{amsart}

\makeatletter
\@namedef{subjclassname@2020}{\textup{2020} Mathematics Subject Classification}
\makeatother

\usepackage{geometry}
\title{On Zappa--Sz\'ep products of multiplicative racks}
\author{Akihiro Takano}

\subjclass[2020]{Primary 20M10; Secondary 20E22, 20M12}
\keywords{multiplicative rack, Zappa--Sz\'ep product, universal group, Ore monoid}

\address{Department of Mathematics, Graduate School of Science
Osaka University,
1-1 Machikaneyama, Toyonaka,
Osaka, 560-0043, Japan}
\email{takano.akihiro.sci@osaka-u.ac.jp}

\usepackage{amsmath} 
\usepackage{amsfonts}
\usepackage{mathrsfs}
\usepackage{amssymb}
\usepackage{amsthm}
\usepackage{mathtools}
\usepackage{bm}
\usepackage[new]{old-arrows}
\usepackage{wrapfig}
\usepackage{graphicx}
\usepackage{color}
\usepackage[labelsep=colon]{caption}
\usepackage[labelformat=empty,subrefformat=parens]{subcaption}
\usepackage{multicol}
\usepackage{units}
\usepackage[all]{xy}
\usepackage{tikz-cd}
\usepackage{hyperref}
\usepackage{mathpazo}
\usepackage{autobreak}
\usepackage{pdfpages}
\usepackage{fancyhdr}
\usepackage{enumitem}

\hypersetup{
colorlinks=true, 
linkcolor=blue, 
citecolor=purple, 
}

\newtheorem{thm}{Theorem}[section]
\newtheorem{prop}[thm]{Proposition}
\newtheorem{lemma}[thm]{Lemma}
\newtheorem{cor}[thm]{Corollary}

\newtheorem*{thm*}{Main Theorem}

\theoremstyle{definition}
\newtheorem{defi}[thm]{Definition}
\newtheorem{ex}[thm]{Example}
\newtheorem{remark}[thm]{Remark}

\newcommand{\bN}{\mathbb{N}}
\newcommand{\bZ}{\mathbb{Z}}

\newcommand{\bF}{\mathbb{F}}

\newcommand{\cM}{\mathcal{M}}

\newcommand{\cG}{\mathcal{G}}
\newcommand{\cB}{\mathcal{B}}

\newcommand{\Map}{\operatorname{Map}}
\newcommand{\Hom}{\operatorname{Hom}}

\newcommand{\Conj}{\operatorname{Conj}}
\newcommand{\Aut}{\operatorname{Aut}}
\newcommand{\Mon}{\operatorname{Mon}}
\newcommand{\Gr}{\operatorname{Gr}}
\newcommand{\id}{\operatorname{id}}
\newcommand{\pr}{\operatorname{pr}}
\newcommand{\Stab}{\operatorname{Stab}}

\newcommand{\can}{\operatorname{can}}
\newcommand{\red}{\operatorname{red}}
\newcommand{\op}{\operatorname{op}}

\newcommand{\wtil}{\widetilde}
\newcommand{\what}{\widehat}

\newcommand{\tleft}{\triangleleft}
\newcommand{\tright}{\triangleright}

\renewcommand{\MR}{\operatorname{MR}}

\newcommand\cinput[2]{\lower#1pt\hbox{\input{#2}}}

\makeatletter
\newcommand{\conditem}[1]{%
  \item[{\rm \textbf{(#1)}}]%
  \phantomsection%
  \protected@edef\@currentlabel{#1}%
}
\makeatother

\usepackage[pagewise]{lineno}

\begin{document}

\begin{abstract}
We study Zappa--Sz\'ep products of multiplicative racks, namely monoids equipped with compatible rack operations.
We characterize multiplicative rack structures on Zappa--Sz\'ep products of monoids and prove that every exact factorization of multiplicative racks induces a canonical exact factorization of their universal groups.
We also establish exact factorizations of groups of units and cancellative quotients, show that reduction modulo units yields a direct product, and obtain consequences for Ore localization.
We then generalize Bianchi's construction using augmentations from partially multiplicative racks into arbitrary multiplicative racks and establish a reconstruction theorem for multiplicative racks with suitable unital ideals.
We investigate the interaction between unital ideals and exact factorizations, and characterize Zappa--Sz\'ep factorizations of generalized Bianchi-type multiplicative racks.
Finally, we use a layered central-quotient construction to obtain finite and infinite non-cancellative multiplicative quandles with non-group factors and nontrivial Zappa--Sz\'ep actions, all directly indecomposable as monoids.
\end{abstract}

\maketitle

\section{Introduction}\label{sec:introduction}
Every group has a natural quandle structure given by conjugation.
Racks and quandles abstract the self-distributive properties of conjugation and have played an important role in knot theory since the work of Joyce~\cite{Joyce1982} and Matveev~\cite{Matveev1982}.
When a multiplication is also present, it is natural to ask how the two operations interact.
Multiplicative racks provide a setting for this question:
they are monoids equipped with rack operations compatible with multiplication, generalizing the conjugation structure of groups.
Since the underlying monoids of multiplicative racks need not be cancellative, it is natural to ask which aspects of their structure are preserved under universal group completion.

Structures combining multiplication and self-distributive operations have appeared in several contexts.
Examples include LD-monoids studied by Dehornoy~\cite{Dehornoy1997}, enveloping monoidal quandles introduced by Kamada and Matsumoto~\cite{KamadaMatsumoto2005}, and qualgebras considered by Lebed~\cite{Lebed2015}.
Bianchi developed partially multiplicative quandles and also considered partially multiplicative racks~\cite{Bianchi2025}.
Multiple conjugation quandles and multiple group racks provide further related structures arising in the study of handlebody-knots and spatial surfaces~\cite{Ishii2015,IshiiMatsuzakiMurao2020}.
In this paper, we focus on multiplicative racks with globally defined multiplication and investigate their Zappa--Sz\'ep factorizations.

An exact factorization $M=A\bowtie U$ of a monoid means that every element of $M$ can be written uniquely as $au$ with $a\in A$ and $u\in U$.
It determines two component actions through the relation
\begin{align*}
ua=(u\tright a)(u\tleft a).
\end{align*}
Zappa--Sz\'ep products generalize semidirect products of groups without requiring either factor to be normal~\cite{Brin2005,Zenab2018}, and have also been studied for monoids, including Garside monoids~\cite{GebhardtTawn2016}.
When the factors are multiplicative racks, a natural problem is to determine when their Zappa--Sz\'ep product admits a compatible rack operation.
We characterize such operations in terms of monoid automorphisms extending the right translations of the factors (Theorem~\ref{thm:external ZSofMR}).
We also show that both component actions are necessarily bijective (Proposition~\ref{prop:MR ZSaction bijective}).

Our main question is whether exact factorizations are preserved under universal group completion.
Let $\cG(M)$ denote the universal group of a monoid $M$, with canonical homomorphism $\omega_M\colon M\to\cG(M)$.
For a multiplicative rack $M$, the map $\omega_M$ preserves the rack operation when $\cG(M)$ is equipped with conjugation, and the universal group construction satisfies the natural adjunction
\begin{align*}
\Hom_{\Gr}(\cG(M),G)\cong\Hom_{\MR}(M,\cM(G)).
\end{align*}
For general monoids, adjoining inverses need not preserve exact factorizations.
The situation is different for multiplicative racks.
Our main result (Theorem~\ref{thm:universal-exactness}) shows that every exact factorization $M=A\bowtie U$ of multiplicative racks induces a canonical exact factorization
\begin{align*}
\cG(M)\cong\cG(A)\bowtie\cG(U).
\end{align*}
The proof uses the kernel congruence of the universal group map and the constraints imposed by the rack operation on the Zappa--Sz\'ep actions.
These constraints also yield an exact factorization of the group of units and a direct-product decomposition after reduction modulo units:
\begin{align*}
M^\times=A^\times\bowtie U^\times,\quad M_{\red}\cong A_{\red}\times U_{\red}.
\end{align*}
We further obtain corresponding results for canonical cancellative quotients and, in the cancellative case, groups of fractions.

A second aim of this paper is to understand non-cancellative multiplicative racks through their ideals.
In particular, minimal ideals of multiplicative racks are groups, but restricting attention to such ideals does not capture the full range of non-cancellative structures.
Motivated by Bianchi's disjoint-union construction~\cite{Bianchi2025}, we introduce a generalized construction based on an augmentation from a partially multiplicative rack into an arbitrary multiplicative rack.
We also prove a reconstruction theorem:
a multiplicative rack with a proper unital ideal invariant under all right translations can be recovered from its complement and the ideal (Theorem~\ref{thm:MR-Bianchi}).
Thus unital ideals provide a framework for studying multiplicative racks beyond the setting of minimal group ideals.

We then investigate how these ideal-theoretic constructions interact with exact factorizations.
We describe the factorization induced on a unital ideal when it meets one of the factors (Proposition~\ref{prop:unital-ideal-factors}).
Using this result, we characterize exact Zappa--Sz\'ep factorizations of generalized Bianchi-type multiplicative racks whose first factor meets the distinguished ideal (Theorem~\ref{thm:Bianchi-ZS-factorization}).
In particular, when the second factor is a group, a suitable exact factorization of the ideal extends to the whole multiplicative rack.
We conclude with a layered construction based on exact factorizations of groups and families of central subgroups indexed by commutative monoids.
This construction yields finite and infinite non-cancellative multiplicative quandles with non-group factors and nontrivial Zappa--Sz\'ep actions, and allows us to compute their universal groups.
All the resulting examples are directly indecomposable as monoids, while the infinite example has a proper unital ideal but no minimal ideal.

The paper is organized as follows.
Section~\ref{sec:preliminaries} recalls the necessary preliminaries.
Section~\ref{sec:definition} characterizes multiplicative rack structures on Zappa--Sz\'ep products.
Section~\ref{sec:universal} studies universal groups, unit groups, and related factorizations.
Section~\ref{sec:non-cancellative} studies unital ideals and generalized Bianchi-type factorizations, and concludes with a layered construction yielding finite and infinite examples.

\subsection*{Acknowledgements}
The author would like to thank Seiichi Kamada and Jumpei Yasuda for their valuable mathematical comments and helpful suggestions.
The author was partially supported by JSPS KAKENHI Grant Number 24KJ0144.

\subsection*{Use of AI}
OpenAI's ChatGPT was used to assist in developing mathematical arguments, drafting portions of the manuscript, and revising the exposition.
The author has reviewed the mathematical content and takes full responsibility for the results and any remaining errors.

\section{Preliminaries} \label{sec:preliminaries}
In this section, we recall racks and monoids, introduce multiplicative racks, and review Zappa--Sz\'ep products of monoids.

\subsection{Racks}  \label{subsec:rack}
We follow the standard conventions for racks and quandles~\cite{Joyce1982,Matveev1982}.

\begin{defi}\label{def:rack}
A \textbf{rack} is a nonempty set $X$ equipped with a binary operation $* \colon X \times X \to X$ satisfying the following conditions:
\begin{enumerate}
\conditem{R1} \label{R1} for all $y \in X$, the \textbf{right translation} $S_y \colon X \to X; x \mapsto x*y$ is bijective,
\conditem{R2} \label{R2} $(x*y)*z = (x*z)*(y*z) \quad (x,y,z \in X)$.
\end{enumerate}
A rack $X$ is a \textbf{quandle} if the condition
\begin{enumerate}
\conditem{Q} \label{Q} $x*x = x \quad (x \in X)$
\end{enumerate}
is also satisfied.
A rack $X$ is \textbf{unital} if there exists $e \in X$, called a \textbf{unit}, such that the following holds:
\begin{enumerate}
\conditem{U} \label{U} $x*e = x$ and $e*x = e \quad (x \in X)$.
\end{enumerate}
\end{defi}

A rack $X$ is \textbf{trivial} if $x*y=x$ for all $x,y\in X$.
Equivalently, the map $S_y$ is the identity for all $y \in X$.
Note that the unit of a unital rack need not be unique.
For example, every element of a trivial rack is a rack unit.

A subset $A$ of $X$ is called a (unital) \textbf{subrack} if $A$ has a (unital) rack structure by restricting the underlying binary operation $*$.
A map $f \colon X \to Y$ between racks is a \textbf{homomorphism} if $f(x*y) = f(x) * f(y)$ for any $x,y \in X$.
A bijective rack homomorphism is an \textbf{isomorphism}.

\begin{ex}\label{ex:conjugation quandle}
Any group $G$ has a unital quandle structure with a binary operation $* \colon G \times G \to G$ defined by
\begin{align*}
g*h \coloneq h^{-1}gh \quad (g,h \in G).
\end{align*}
This is the \textbf{conjugation quandle} of $G$, denoted by $\Conj(G)$.
\end{ex}

\subsection{Monoids} \label{subsec:monoid}
\begin{defi}\label{def:monoid}
A \textbf{monoid} is a set $M$ equipped with an associative binary operation $\cdot \colon M \times M \to M$ and a two-sided identity element $1_M \in M$.
That is, for all $x,y,z \in M$,
\begin{align*}
(x \cdot y) \cdot z = x \cdot (y \cdot z) \quad\text{and}\quad
x \cdot 1_M = 1_M \cdot x = x.
\end{align*}
\end{defi}

Throughout this paper, we omit the multiplication symbol $\cdot$ and write $xy$ instead of $x \cdot y$.
We write $1_M$ for the identity of a monoid $M$, or simply $1$ when no confusion can arise.
An element $u \in M$ is called a \textbf{unit} if there exists $v \in M$ such that $uv = vu = 1$.
Then $v$ is the inverse of $u$, denoted by $u^{-1}$.
The set of all units of $M$ forms a group, denoted by $M^\times$.
A map $f\colon M\to N$ between monoids is a \textbf{homomorphism} if $f(xy)=f(x)f(y)$ for all $x,y\in M$ and $f(1_M)=1_N$.
A bijective monoid homomorphism is an \textbf{isomorphism}.
For a monoid $M$, its \textbf{opposite monoid} $M^{\op}$ has the same underlying set as $M$, with multiplication defined by $x\cdot_{\op}y\coloneq yx$.

A monoid $M$ is said to be \textbf{left cancellative} if $xy = xz$ implies $y = z$ for all $x, y, z \in M$.
Two elements $x, y \in M$ have a \textbf{common left multiple} $m \in M$ if there exist $z, w \in M$ such that $zx = wy = m$.
The notions of \textbf{right cancellative} and \textbf{common right multiple} are defined similarly.
A monoid $M$ is said to be \textbf{cancellative} if it is both left and right cancellative.
Note that if the elements $u$ and $v$ of a cancellative monoid $M$ satisfy $uv = 1$, then they are units in $M$.
We say that $M$ has \textbf{common left} (resp. \textbf{right}) \textbf{multiples} if every pair of elements has a common left (resp. right) multiple.

For any monoid $M$, there exists a group $\cG(M)$ and a monoid homomorphism $\omega_M \colon M \to \cG(M)$ with the following universal property:
for any group $H$ and any monoid homomorphism $\varphi \colon M \to H$, there exists a unique group homomorphism $\wtil{\varphi} \colon \cG(M) \to H$ such that $\varphi = \wtil{\varphi} \circ \omega_M$.
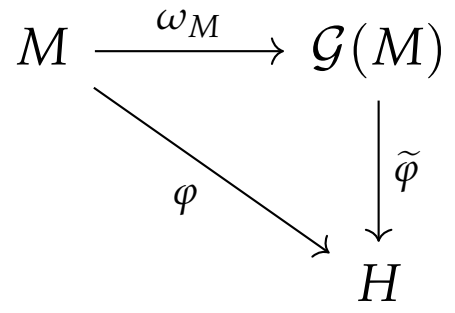
\begin{center}
\begin{tikzcd}[sep = large]
M \arrow[r, "\omega_M"] \arrow[rd, "\varphi"'] & \cG(M) \arrow[d, "\wtil{\varphi}"]\\
& H
\end{tikzcd}
\end{center}
This group $\cG(M)$ is unique up to group isomorphism, and is called the \textbf{universal group} of $M$.
Explicitly, $\cG(M)$ is generated by the elements of $M$, subject to the relations defining its multiplication.
The map $\omega_M$ is injective if and only if $M$ embeds into a group.
For Ore monoids, the universal group admits a description in terms of fractions.

\begin{defi}\label{def:Ore}
An \textbf{Ore monoid} is a cancellative monoid with common right multiples.
\end{defi}

By Ore's theorem (see \cite[Theorem~1.23]{CliffordPreston1961}), an Ore monoid $M$ embeds in its \textbf{group of right fractions} $\what{M}$, in which every element is of the form $mn^{-1}$ with $m,n\in M$.
The universal property of localization gives a canonical isomorphism $\cG(M)\cong\what{M}$.

\subsection{Multiplicative racks}
Multiplicative racks are precisely the complete partially multiplicative racks considered by Bianchi~\cite[Appendix]{Bianchi2025}.
The compatibility identities of group multiplication and conjugation motivate the following definition.

\begin{defi}\label{def:MR}
A \textbf{multiplicative rack} is a monoid $M$ equipped with a rack operation $* \colon M \times M \to M$ such that the following conditions hold:
\begin{enumerate}
\conditem{MR1}\label{MR1} $xy = y(x*y) \quad (x,y \in M)$,
\conditem{MR2}\label{MR2} $x*(yz) = (x*y)*z \quad (x,y,z \in M)$,
\conditem{MR3}\label{MR3} $(xy)*z = (x*z)(y*z) \quad (x,y,z \in M)$.
\end{enumerate}
A multiplicative rack that is also a quandle is called a \textbf{multiplicative quandle}.
\end{defi}

Let $M$ be a multiplicative rack.
Since $S_1=S_1 \circ S_1$ by condition~\eqref{MR2} and $S_1$ is bijective, we have $S_1=\id_M$.
By condition~\eqref{MR3}, every right translation is a monoid automorphism, and thus $S_x(1) = 1$ for any $x \in M$.
Therefore, $1$ is a rack unit of $M$.

A \textbf{multiplicative subrack} is a subset that is both a submonoid and a subrack, with the induced operations and common identity.
A map $f\colon M\to N$ between multiplicative racks is a \textbf{homomorphism} if it is both a monoid homomorphism and a rack homomorphism.
A bijective homomorphism is called an \textbf{isomorphism}.

\begin{remark}\label{remark:MR}
A monoid $M$ admits a multiplicative rack structure if and only if there exists a monoid anti-homomorphism
\begin{align*}
S \colon M \to \Aut_{\Mon}(M)
\end{align*}
such that $xy=yS_y(x)$ for all $x,y\in M$.
The corresponding rack operation is given by $x*y\coloneq S_y(x)$.
Indeed, the necessity follows from conditions~\eqref{MR1}--\eqref{MR3} and the bijectivity of right translations.
Conversely, such an anti-homomorphism defines a binary operation satisfying conditions~\eqref{MR1}--\eqref{MR3}, whose right translations are bijective.
Moreover, since $yz=z(y*z)$, we have
\begin{align*}
S_z\circ S_y=S_{yz}=S_{z(y*z)}=S_{y*z}\circ S_z,
\end{align*}
which gives condition~\eqref{R2}.
\end{remark}

\begin{ex}\label{ex:group MR}
For a group $G$, let $\cM(G)$ denote the \textbf{conjugation multiplicative quandle} of $G$, whose multiplication is the group multiplication and whose rack operation is given by conjugation.
\end{ex}

\begin{ex}\label{ex:comm MR}
Every commutative monoid $M$ admits a multiplicative quandle structure given by the trivial quandle operation.
\end{ex}

\begin{remark} \label{rem:LD}
Multiplicative racks are closely related to Dehornoy's LD-monoids~\cite[Definition~2.1]{Dehornoy1997}.
Indeed, setting $x\circ y=yx$ and $x \star y=y*x$, the operations $\circ$ and $\star$ satisfy the LD-monoid identities.
In our setting, the left translations of $\star$ are additionally bijective.
\end{remark}

\subsection{Zappa--Sz\'ep products} \label{subsec:ZS}
We recall the standard monoid Zappa--Sz\'ep product; see Brin~\cite{Brin2005} and Zenab~\cite{Zenab2018}.

\begin{defi}\label{def:ZS}
Let $M$ be a monoid, and let $A$ and $U$ be submonoids of $M$.
We say that $M$ is an \textbf{internal Zappa--Sz\'ep product} of $A$ and $U$, or equivalently, that $M=AU$ is an \textbf{exact factorization}, if the multiplication map
\begin{align*}
A \times U \to M;\ (a,u) \to au
\end{align*}
is bijective.
In this case, we write $M = A \bowtie U$.
Note that if $M$ is a group and $A$ and $U$ are its subgroups, then this is equivalent to $M = AU$ and $A \cap U = \{ 1_M \}$.
\end{defi}

An internal exact factorization identifies $M$ with $A\times U$ as a set.
For $a\in A$ and $u\in U$, there are unique $b\in A$ and $v\in U$ such that $ua=bv$.
We write $u\tright a=b$ and $u\tleft a=v$.
Then we obtain two maps
\begin{align*}
\begin{aligned}
\tright &\colon U\times A\to A;\ (u,a)\mapsto u\tright a,\\
\tleft &\colon U\times A\to U;\ (u,a)\mapsto u\tleft a.
\end{aligned}
\end{align*}

\begin{lemma}\label{lemma:ZS}
Let $M = A \bowtie U$ be an internal Zappa--Sz\'ep product of submonoids $A$ and $U$.
Then the maps
\begin{align*}
\tright \colon U \times A \to A \quad\text{and}\quad \tleft \colon U \times A \to U
\end{align*}
defined above satisfy the following conditions$:$
\begin{enumerate}
\conditem{ZS1}\label{ZS1} $1_U  \tright a = a \quad (a \in A)$,
\conditem{ZS2}\label{ZS2} $1_U \tleft a = 1_U \quad (a \in A)$,
\conditem{ZS3}\label{ZS3} $u \tleft 1_A = u \quad (u \in U)$,
\conditem{ZS4}\label{ZS4} $u \tright 1_A = 1_A \quad (u \in U)$,
\conditem{ZS5}\label{ZS5} $(uv) \tright a = u \tright (v \tright a) \quad (a \in A, u,v \in U)$,
\conditem{ZS6}\label{ZS6} $u \tright (ab) = (u \tright a)((u \tleft a) \tright b) \quad (a,b \in A, u \in U)$,
\conditem{ZS7}\label{ZS7} $u \tleft (ab) = (u \tleft a) \tleft b \quad (a,b \in A, u \in U)$,
\conditem{ZS8}\label{ZS8} $(uv) \tleft a = (u \tleft (v \tright a))(v \tleft a) \quad (a \in A, u,v \in U)$.
\end{enumerate}
\end{lemma}

\begin{proof}
Conditions~\eqref{ZS1} and \eqref{ZS2} follow from $a = 1_U a = (1_U \tright a) (1_U \tleft a)$.
Similarly, conditions~\eqref{ZS3} and \eqref{ZS4} follow from $u = u1_A = (u \tright 1_A)(u \tleft 1_A)$.
Moreover,
\begin{align*}
(uv)a = u(va) = u(v \tright a)(v \tleft a) = (u \tright (v \tright a))(u \tleft (v \tright a))(v \tleft a)
\end{align*}
and
\begin{align*}
u(ab) = (ua)b = (u \tright a)(u \tleft a)b = (u \tright a)((u \tleft a) \tright b)((u \tleft a) \tleft b)
\end{align*}
yield conditions~\eqref{ZS5}--\eqref{ZS8} by comparing the two components of each factorization.
\end{proof}

Under this identification, the multiplication of $M$ is given on $A\times U$ by
\begin{align*}
(a,u)(b,v)=(a(u\tright b),(u\tleft b)v).
\end{align*}
Thus the map $(a,u)\mapsto au$ is a monoid isomorphism with respect to this multiplication.

\begin{ex}\label{ex:semidirect}
Let $G=AU$ be an exact factorization of groups with $U\trianglelefteq G$.
Then $ua=a(a^{-1}ua)$ for $a\in A$ and $u\in U$, so $u\tright a=a$ and $u\tleft a=a^{-1}ua$.
\end{ex}

The preceding lemma shows that every internal Zappa--Sz\'ep product determines a pair of maps satisfying conditions~\eqref{ZS1}--\eqref{ZS8}.
Conversely, these conditions allow us to construct a monoid from two given monoids.

\begin{defi} \label{def:external ZS}
Let $A$ and $U$ be monoids.
A pair of maps
\begin{align*}
\tright &\colon U\times A\to A;\ (u,a)\mapsto u\tright a,\\
\tleft &\colon U\times A\to U;\ (u,a)\mapsto u\tleft a
\end{align*}
is called a \textbf{Zappa--Sz\'ep action} if the two maps satisfy conditions~\eqref{ZS1}--\eqref{ZS8}.
Then the set $A \times U$ equipped with the multiplication defined by
\begin{align*}
(a,u) (b,v) \coloneq (a(u \tright b), (u \tleft b)v)
\end{align*}
is called the \textbf{external Zappa--Sz\'ep product} of $A$ and $U$, denoted by $A \bowtie U$.
\end{defi}

Conditions~\eqref{ZS1}--\eqref{ZS8} ensure that $A\bowtie U$ is a monoid with identity $(1_A,1_U)$.
Moreover, the canonical embeddings
\begin{align*}
A\hookrightarrow A\bowtie U;\ a\mapsto(a,1_U)
\quad\text{and}\quad
U\hookrightarrow A\bowtie U;\ u\mapsto(1_A,u)
\end{align*}
identify the external product with an internal exact factorization.

Let $A$ and $U$ be monoids with a Zappa--Sz\'ep action.
For each $a \in A$ and $u \in U$, we can define maps
\begin{align*}
R_a \colon U \to U;\ v \mapsto v \tleft a \quad\text{and}\quad L_u \colon A \to A;\ b \mapsto u \tright b.
\end{align*}
In general, these maps are not monoid homomorphisms, as conditions~\eqref{ZS6} and \eqref{ZS8} show.
For a set $X$, let $\Map(X)$ denote the monoid of all maps from $X$ to itself under composition.
Conditions~\eqref{ZS1}, \eqref{ZS3}, \eqref{ZS5}, and \eqref{ZS7} imply that
\begin{align*}
L\colon U\to\Map(A);\ u\mapsto L_u
\end{align*}
is a monoid homomorphism, whereas
\begin{align*}
R\colon A\to\Map(U);\ a\mapsto R_a
\end{align*}
is a monoid anti-homomorphism.
For arbitrary Zappa--Sz\'ep products of monoids, neither $R_a$ nor $L_u$ need be injective or surjective.

\begin{remark}\label{remark:group ZS}
Let $M = A \bowtie U$ be an exact factorization of monoids.
If $M$ is a group, then both $A$ and $U$ are also groups.
Conversely, if $A$ and $U$ are groups, then the external Zappa--Sz\'ep product $A \bowtie U$ is also a group.
The inverse of $(a,u)\in A\bowtie U$ is $(u^{-1} \tright a^{-1}, u^{-1} \tleft a^{-1})$.
Moreover, the maps $R_a \colon U \to U$ and $L_u \colon A \to A$ are both bijective for any $a \in A$ and $u \in U$.
\end{remark}

\section{Zappa--Sz\'ep products of multiplicative racks} \label{sec:definition}
In this section, we characterize multiplicative rack structures on Zappa--Sz\'ep products of monoids that extend the given structures on the factors.

\begin{defi}\label{def:ZSofMR}
Let $M$ be a multiplicative rack and $A$ and $U$ be multiplicative subracks of $M$.
We call $M=A\bowtie U$ an \textbf{internal Zappa--Sz\'ep product} of multiplicative racks if the underlying monoid is an internal Zappa--Sz\'ep product of $A$ and $U$.
\end{defi}

By Remark~\ref{remark:MR}, the right translations of $M$ define a monoid anti-homomorphism
\begin{align*}
S\colon M\to\Aut_{\Mon}(M)
\end{align*}
satisfying $pq=qS_q(p)$ for all $p,q\in M$.
Since every $p\in M$ has a unique expression $p=au$ with $a\in A$ and $u\in U$, we have $S_p=S_u\circ S_a$.

For $a\in A$ and $u\in U$, put $\Phi_a\coloneq S_a$ and $\Psi_u\coloneq S_u$.
Then
\begin{align*}
\Phi\colon A\to\Aut_{\Mon}(M)
\quad\text{and}\quad
\Psi\colon U\to\Aut_{\Mon}(M)
\end{align*}
are monoid anti-homomorphisms.
Moreover, the Zappa--Sz\'ep relation $ua=(u\tright a)(u\tleft a)$ gives
\begin{align*}
\Phi_a\circ\Psi_u
=
\Psi_{u\tleft a}\circ\Phi_{u\tright a}
\end{align*}
for all $a\in A$ and $u\in U$.
Under the identification of $M$ with $A\times U$ as sets, the multiplication and rack operation are given by
\begin{align*}
(a,u)\cdot(b,v)
&=(a(u\tright b),(u\tleft b)v),\\
(a,u)*(b,v)
&=(\Psi_v\circ\Phi_b)(a,u)
\end{align*}
for all $a,b\in A$ and $u,v\in U$.
Thus the map $A\times U\to M;\ (a,u)\mapsto au$, is an isomorphism of multiplicative racks with respect to these operations.

Conversely, given a Zappa--Sz\'ep product of monoids, we ask when it admits a multiplicative rack structure extending those of the factors.
The following theorem characterizes such structures in terms of monoid automorphisms.

\begin{thm}\label{thm:external ZSofMR}
Let $A$ and $U$ be multiplicative racks, and let $M=A\bowtie U$ be the external Zappa--Sz\'ep product of their underlying monoids.
Assume that the underlying Zappa--Sz\'ep action
\begin{align*}
\tright \colon U \times A \to A \quad\text{and}\quad \tleft \colon U \times A \to U
\end{align*}
is fixed.
The following conditions are equivalent$:$
\begin{enumerate}
\item There exists a multiplicative rack structure on $M$ such that the canonical inclusions $A \hookrightarrow A \times U$ and $U \hookrightarrow A \times U$ are homomorphisms of multiplicative racks.
\item There exist monoid anti-homomorphisms
\begin{align*}
\Phi \colon A \to \Aut_{\Mon}(M) \quad\text{and}\quad \Psi \colon U \to \Aut_{\Mon}(M)
\end{align*}
satisfying the following conditions:
\begin{enumerate}
\conditem{E1}\label{E1} $\Phi_b(a,1_U) = (a*b, 1_U)$ and $\Psi_v(1_A,u) = (1_A, u*v) \quad (a,b \in A, u,v \in U)$,
\conditem{E2}\label{E2} $\Phi_a \circ \Psi_u = \Psi_{u \tleft a} \circ \Phi_{u \tright a} \quad (a \in A, u \in U)$,
\conditem{E3}\label{E3} $p(a,1_U) = (a,1_U)\Phi_a(p)$ and $p(1_A,u) = (1_A,u)\Psi_u(p) \quad (p \in M, a \in A, u \in U)$.
\end{enumerate}
\end{enumerate}
For each such pair $(\Phi,\Psi)$, the multiplicative rack operation is uniquely determined by
\begin{align*}
(a,u) * (b,v) \coloneq (\Psi_v \circ \Phi_b)(a,u).
\end{align*}
When these equivalent conditions hold, we call $M$ the \textbf{external Zappa--Sz\'ep product} of multiplicative racks, and denote it again by $M = A \bowtie U$.

Moreover, $M=A\bowtie U$ is a multiplicative quandle if and only if the anti-homomorphisms $\Phi$ and $\Psi$ satisfy the additional condition
\begin{enumerate}
\conditem{EQ}\label{EQ} $(\Psi_u \circ \Phi_a)(a,u) = (a,u) \quad (a \in A, u \in U)$.
\end{enumerate}
\end{thm}

\begin{proof}
\underline{(1) $\Rightarrow$ (2)}\ \ Under (1), $M$ is an internal Zappa--Sz\'ep product of $A \times \{1_U\}$ and $\{1_A\} \times U$.
Define $\Phi_a=S_{(a,1_U)}$ and $\Psi_u=S_{(1_A,u)}$.
These maps satisfy conditions~\eqref{E1}--\eqref{E3} by the preceding discussion.\\
\underline{(2) $\Rightarrow$ (1)}\ \
For $q=(b,v) \in M = A \times U$, we define $S_q \coloneq \Psi_v \circ \Phi_b$ and $p*q \coloneq S_q(p)$.
We verify the conditions in Remark~\ref{remark:MR}.
For $p=(a,u)$ and $q=(b,v)$, we have $pq=(a(u\tright b),(u\tleft b)v)$ and hence
\begin{align*}
S_{pq} &= \Psi_{(u \tleft b)v} \circ \Phi_{a(u \tright b)}\\
&= \Psi_v \circ \Psi_{u \tleft b} \circ \Phi_{u \tright b} \circ \Phi_a\\
&= \Psi_v \circ \Phi_b \circ \Psi_u \circ \Phi_a = S_q \circ S_p,
\end{align*}
where we used that $\Phi$ and $\Psi$ are anti-homomorphisms and condition~\eqref{E2}.
Since each $S_q$ is a monoid automorphism, Remark~\ref{remark:MR} shows that it remains only to verify $pq=qS_q(p)$.
Write $q=(b,v)=(b,1_U)(1_A,v)$.
Condition~\eqref{E3} yields
\begin{align*}
pq=(b,1_U)\Phi_b(p)(1_A,v)=(b,1_U)(1_A,v)\Psi_v(\Phi_b(p))=qS_q(p).
\end{align*}
Consequently, $p*q=S_q(p)$ is a multiplicative rack operation.
By condition~\eqref{E1}, we have
\begin{align*}
(a,1_U)*(b,1_U)=\Phi_b(a,1_U)=(a*b,1_U),
\end{align*}
so the canonical inclusion $A \hookrightarrow A \times U$ is a rack homomorphism.
Similarly, $U \hookrightarrow A \times U$ is also a rack homomorphism.

Finally, $M$ is a multiplicative quandle if and only if $q*q=q$ for every $q=(a,u)\in M$, which is precisely condition~\eqref{EQ}.
\end{proof}

The multiplicative rack structure forces both component maps of the Zappa--Sz\'ep action to be bijective.

\begin{prop}\label{prop:MR ZSaction bijective}
Let $M = A \bowtie U$ be a Zappa--Sz\'ep product of multiplicative racks.
Then the maps
\begin{align*}
R_a \colon U \to U;\ u \mapsto u \tleft a \quad\text{and}\quad L_u \colon A \to A;\ a \mapsto u \tright a
\end{align*}
are bijective for all $a \in A$ and $u \in U$.
\end{prop}

\begin{proof}
Let $\Phi_a=S_{(a,1_U)}$ and $\Psi_u=S_{(1_A,u)}$ be the monoid automorphisms from Theorem~\ref{thm:external ZSofMR}.\\
\underline{Bijectivity of $R_a$}:
Fix $a \in A$, and define the maps $\eta_a, \eta'_a \colon U \to U$ by
\begin{align*}
\eta_a(u) \coloneq \pr_U(\Phi_a(1_A,u)),\quad \eta'_a(u) \coloneq \pr_U(\Phi_a^{-1}(1_A,u))
\end{align*}
Writing $\Phi_a(1_A,u) = (b, \eta_a(u))$ for some $b \in A$, we obtain
\begin{align*}
(1_A, u) = \Phi_a^{-1} (\Phi_a(1_A, u)) = \Phi_a^{-1}(b, 1_U)\Phi_a^{-1}(1_A, \eta_a(u)).
\end{align*}
By condition~\eqref{E1}, both $\Phi_a$ and $\Phi_a^{-1}$ preserve $A \times \{ 1_U \}$.
Comparing the $U$-coordinates gives $u = \eta'_a(\eta_a(u))$, and the reverse composition is treated similarly.
Thus $\eta_a$ and $\eta'_a$ are mutually inverse.
From condition~\eqref{E3},
\begin{align*}
(u\tright a,u\tleft a) =(1_A,u)(a,1_U) =(a,1_U)\Phi_a(1_A,u),
\end{align*}
so $R_a=\eta_a$ is bijective.\\
\underline{Bijectivity of $L_u$}:
Similarly, fix $u \in U$, and define the maps $\theta_u, \theta'_u \colon A \to A$ by
\begin{align*}
\theta_u(a) \coloneq \pr_A(\Psi_u(a,1_U)),\quad \theta'_u(a) \coloneq \pr_A(\Psi_u^{-1}(a,1_U))
\end{align*}
Since both $\Psi_u$ and $\Psi_u^{-1}$ preserve $\{ 1_A \} \times U$ by condition~\eqref{E1}, the same argument shows that $\theta_u$ and $\theta_u'$ are mutually inverse.
From condition~\eqref{E3},
\begin{align*}
(a,u) = (a,1_U)(1_A,u) = (1_A,u)\Psi_u(a,1_U),
\end{align*}
which gives $a=u\tright\theta_u(a)=L_u(\theta_u(a))$.
Thus $L_u (= \theta'_u)$ is bijective.
\end{proof}

\section{Universal groups and exact factorizations}\label{sec:universal}
In this section, we show that exact factorizations of multiplicative racks induce exact factorizations of their universal groups.
We first establish the necessary structural properties and describe the kernel congruence of the canonical map, then derive consequences for unit groups, cancellative quotients, and Ore localization.

\subsection{Structural properties of multiplicative racks}\label{subsec:properties-MR}
\begin{lemma} \label{lemma:MRidempotent}
Let $M$ be a multiplicative rack.
Then the following hold$:$
\begin{enumerate}
\item $Mx = xM$ for every $x \in M$,
\item Every idempotent $e \in M$ is central, and satisfies $x*e = x$ for all $x \in M$.
\end{enumerate}
\end{lemma}

\begin{proof}
(1)\ By condition~\eqref{MR1} and the surjectivity of $S_x$, we have $Mx = x(S_x(M)) = xM$.\\
(2)\ Suppose that $e^2 = e$.
Condition~\eqref{MR2} gives $S_e = S_e \circ S_e$.
Since $S_e$ is bijective, it follows that $S_e = \id_M$.
Consequently, $xe = e(x*e) = ex$, so $e$ is central.
\end{proof}

\begin{lemma}\label{lem:onesided-MR-units}
If $M$ is a multiplicative rack and $ab=1$ in $M$, then $ba=1$.
\end{lemma}

\begin{proof}
Put $e=ba$.
Then $e^2=e$, so Lemma~\ref{lemma:MRidempotent} implies that $e$ is central.
Since $ae=(ab)a=a$, centrality gives $ea=a$ and hence $e=e(ab)=(ea)b=ab=1$.
\end{proof}

\begin{lemma}\label{lemma:MRhasCM}
The underlying monoid of every multiplicative rack has both common right and left multiples.
\end{lemma}

\begin{proof}
For any $x,y\in M$, condition~\eqref{MR1} gives $xy=y(x*y)\in xM\cap yM$, so $M$ has common right multiples.
By Lemma~\ref{lemma:MRidempotent}(1), we also have $xy\in Mx\cap My$.
Thus $M$ has common left multiples.
\end{proof}

\begin{prop}\label{prop:MR+cancellative MQ}
Let $M$ be a multiplicative rack.
Then its underlying monoid is left-cancellative if and only if it is right-cancellative.
Under either condition, $M$ is a multiplicative quandle.
\end{prop}

\begin{proof}
Suppose that $M$ is left-cancellative.
Let $x,y,a \in M$ with $xa=ya$.
From condition~\eqref{MR1}, we have $a(x*a) = a(y*a)$ and thus $x*a=y*a$.
Since the right translation $S_a$ is injective, $x=y$ holds.
Conversely, suppose that $M$ is right-cancellative.
Let $x,y,a\in M$ satisfy $ax=ay$.
By condition~\eqref{MR1}, we have
\begin{align*}
S_a^{-1}(x)a
=a(S_a^{-1}(x)*a)
=ax.
\end{align*}
Similarly, $S_a^{-1}(y)a=ay$.
Right-cancellativity gives $S_a^{-1}(x)=S_a^{-1}(y)$, and hence $x=y$.

If $M$ is left-cancellative, then we have $x^2 = x(x*x)$ for all $x \in M$ from condition~\eqref{MR1}, and thus $x=x*x$ holds.
\end{proof}

In particular, the underlying monoid of any cancellative multiplicative rack is an Ore monoid.
The converse of the quandle assertion fails: a commutative monoid with zero, equipped with the trivial quandle operation, is a non-cancellative multiplicative quandle.
Every finite cancellative monoid is a group.
Consequently, every finite cancellative multiplicative rack is the conjugation multiplicative quandle of a finite group.

\begin{lemma}\label{lemma:cancellative-rigidity}
Let $f \colon M \to N$ be a monoid homomorphism between multiplicative racks.
If $N$ is left-cancellative, then $f$ is automatically a homomorphism of multiplicative racks.
\end{lemma}

\begin{proof}
Applying condition~\eqref{MR1}, we obtain
\begin{align*}
f(q)f(p*q) = f(p)f(q) = f(q)(f(p)*f(q))
\end{align*}
for $p,q \in M$.
Since $N$ is left-cancellative, $f(p*q) = f(p)*f(q)$ holds.
\end{proof}

\begin{remark}
By Lemma~\ref{lemma:cancellative-rigidity}, every monoid homomorphism from a multiplicative rack $M$ to a group $G$ is automatically a homomorphism of multiplicative racks into $\cM(G)$.
In particular, the canonical map $\omega_M\colon M\to\cM(\cG(M))$ preserves both operations.
Hence there is a natural bijection
\begin{align*}
\Hom_{\Gr}(\cG(M),G)
\cong
\Hom_{\MR}(M,\cM(G)).
\end{align*}
Thus $\cG$ is left adjoint to $\cM$.
This is the standard group-completion adjunction in the setting of multiplicative racks; compare Bianchi's enveloping-group construction for partially multiplicative quandles and racks~\cite[Definition~2.9 and Appendix]{Bianchi2025}.
The universal group $\cG(M)$ is also related to the associated group of the underlying rack of $M$, defined by
\begin{align*}
\operatorname{As}(M)
\coloneq
\langle e_x\ (x\in M)\mid e_{x*y}=e_y^{-1}e_xe_y\ (x,y\in M)\rangle.
\end{align*}
By condition~\eqref{MR1}, there is a canonical surjective group homomorphism
\begin{align*}
\operatorname{As}(M)\to\cG(M);\ e_x\mapsto\omega_M(x).
\end{align*}
More precisely, we have a canonical isomorphism
\begin{align*}
\cG(M)\cong
\operatorname{As}(M)\big/
\langle\!\langle e_xe_ye_{xy}^{-1}\mid x,y\in M\rangle\!\rangle,
\end{align*}
where $\langle\!\langle\cdot\rangle\!\rangle$ denotes the normal closure.
\end{remark}

\subsection{The kernel of the universal group map}\label{subsec:universal-groups}
We now study the kernel congruence of the canonical homomorphism $\omega_M$, defined by
\begin{align*}
x \sim_M y \iff \omega_M(x)=\omega_M(y).
\end{align*}
The image $\omega_M(M)$ is a multiplicative subrack of $\cM(\cG(M))$, and is isomorphic to $M_{\can}\coloneq M/{\sim_M}$.
In particular, $M_{\can}$ is cancellative and hence an Ore monoid by Lemma~\ref{lemma:MRhasCM}.
Moreover, every monoid homomorphism from $M$ to a group factors uniquely through the quotient map $M\to M_{\can}$.
Hence there are canonical isomorphisms $\cG(M)\cong\cG(M_{\can})\cong\what{M_{\can}}$.
Since $M_{\can}$ has both common right and left multiples by Lemma~\ref{lemma:MRhasCM}, every element of $\cG(M)$ admits both right- and left-fraction representations.
Thus
\begin{align*}
\cG(M)
&=\{\omega_M(x)\omega_M(y)^{-1}\mid x,y\in M\}\\
&=\{\omega_M(x)^{-1}\omega_M(y)\mid x,y\in M\}.
\end{align*}

\begin{lemma} \label{lem:cancellative-quotient}
Let $M$ be a multiplicative rack.
For $x,y\in M$, the following conditions are equivalent$:$
\begin{enumerate}
\item $x \sim_M y$,
\item There exists $z \in M$ such that $xz = yz$,
\item There exists $w \in M$ such that $wx = wy$.
\end{enumerate}
\end{lemma}

\begin{proof}
We first prove (1) $\Leftrightarrow$ (2).
Define a relation $\approx$ on $M$ by
\begin{align*}
x\approx y \iff xz=yz \text{ for some } z\in M.
\end{align*}
Reflexivity and symmetry are immediate.
Suppose that $xz=yz$ and $yw=tw$.
Since $M$ has common right multiples by Lemma~\ref{lemma:MRhasCM}, there exist $p,q\in M$ such that $zp=wq$.
Then
\begin{align*}
x(zp)=y(zp)=y(wq)=t(wq),
\end{align*}
so $x\approx t$ holds.
Thus $\approx$ is an equivalence relation.
Put $Q\coloneq M/{\approx}$, and let $[x]$ denote the equivalence class of $x\in M$.
We next show that the operations
\begin{align*}
[x][y]\coloneq[xy],\quad [x]*[y]\coloneq[x*y]
\end{align*}
are well-defined and make $Q$ a multiplicative rack.
If $xz=yz$, then $(ax)z=(ay)z$ for every $a\in M$, so $ax \approx ay$.
By Lemma~\ref{lemma:MRhasCM}, we choose $p,q \in M$ such that $zp = aq$.
Then
\begin{align*}
(xa)q=x(aq)=x(zp)=y(zp)=y(aq)=(ya)q.
\end{align*}
Hence $xa\approx ya$, and $\approx$ is a monoid congruence.
Since each right translation $S_a$ is a monoid automorphism, both $S_a$ and $S_a^{-1}$ preserve $\approx$.
Moreover, if $a\approx b$, choose $z\in M$ such that $az=bz$.
By condition~\eqref{MR2}, we have
\begin{align*}
S_z \circ S_a = S_{az} = S_{bz} = S_z \circ S_b.
\end{align*}
Since $S_z$ is bijective, it follows that $S_a=S_b$.
Consequently, the rack operation descends to $Q$, and its right translations remain bijective.
Thus $Q$ is a multiplicative rack.
Moreover, $Q$ is right cancellative.
Indeed, if $[x][z] = [y][z]$, then $xzt = yzt$ for some $t\in M$, and hence $x \approx y$.
By Proposition~\ref{prop:MR+cancellative MQ}, $Q$ is cancellative.
By Lemma~\ref{lemma:MRhasCM}, it has common right multiples and hence is an Ore monoid.
Every monoid homomorphism from $M$ to a group identifies elements related by $\approx$, and hence factors uniquely through the quotient map $\pi \colon M \to Q$.
Therefore, the quotient map $\pi$ induces a canonical
isomorphism $\cG(M) \cong \cG(Q)$ which sends $\omega_M(x)$ to $\omega_Q([x])$ for every $x \in M$.
Since $Q$ is an Ore monoid, its canonical map $\omega_Q \colon Q\to\cG(Q)$ is injective.
Consequently,
\begin{align*}
\omega_M(x)=\omega_M(y) &\iff \omega_Q([x]) = \omega_Q([y])\\
&\iff [x] = [y]\ \text{in $Q$}\\
&\iff x\approx y.
\end{align*}
This proves the equivalence of (1) and (2).

Finally, we prove (2) $\Leftrightarrow$ (3).
If $xz = yz$, then condition~\eqref{MR1} yields $z(x*z) = z(y*z)$.
Applying the monoid automorphism $S_z^{-1}$ gives
\begin{align*}
S_z^{-1}(z)x=S_z^{-1}(z)y.
\end{align*}
Conversely, suppose that $wx=wy$.
By condition~\eqref{MR2}, we have
\begin{align*}
S_x\circ S_w=S_{wx}=S_{wy}=S_y\circ S_w.
\end{align*}
Since $S_w$ is bijective, it follows that $S_x=S_y$.
Therefore, condition~\eqref{MR1} gives
\begin{align*}
xS_x(w)=wx=wy=yS_y(w)=yS_x(w),
\end{align*}
which proves (2).
\end{proof}

\subsection{Exact factorizations of universal groups}\label{subsec:universal-exactness}
We now prove that every exact factorization of multiplicative racks induces an exact factorization of their universal groups.
The proof proceeds in three steps:
injectivity of the canonical homomorphisms, existence of an ordered product decomposition, and uniqueness of the decomposition.

\begin{prop}\label{prop:universal-inclusions}
Let $M=A\bowtie U$ be a Zappa--Sz\'ep product of multiplicative racks.
Then the canonical homomorphisms
\begin{align*}
\cG(A)\to\cG(M),\quad \cG(U)\to\cG(M)
\end{align*}
are injective.
\end{prop}

\begin{proof}
Suppose that $a,b\in A$ have the same image in $\cG(M)$.
By Lemma~\ref{lem:cancellative-quotient} (2), there exists $m\in M$ such that $am=bm$.
Writing $m=cu$ with $c\in A$ and $u\in U$, we obtain $(ac)u=(bc)u$.
Uniqueness of the factorization gives $ac=bc$, so $a$ and $b$ have the same image in $\cG(A)$.
Since every element of $\cG(A)$ can be written as a right fraction $\omega_A(a)\omega_A(b)^{-1}$, the induced homomorphism is injective.

The injectivity of the second homomorphism follows similarly, using Lemma~\ref{lem:cancellative-quotient} (3) and the left-fraction representation of $\cG(U)$.
\end{proof}

We identify $\cG(A)$ and $\cG(U)$ with their images in $\cG(M)$.
The following lemma gives a restriction on the Zappa--Sz\'ep action that will be used in the proof of the main theorem.

\begin{lemma}\label{lem:unit-shifts}
Let $M=A\bowtie U$ be a Zappa--Sz\'ep product of multiplicative racks.
Then
\begin{align*}
u\tright a\in aA^\times
\quad\text{and}\quad
u\tleft a\in U^\times u
\end{align*}
for all $a\in A$ and $u\in U$.
\end{lemma}

\begin{proof}
Fix $a\in A$ and $u\in U$, and write $S_a(u)=pq$ with $p\in A$ and $q\in U$.
By the argument in the proof of Proposition~\ref{prop:MR ZSaction bijective}, there exists $r\in A$ such that $S_a^{-1}(q)=ru$.
Since $S_a^{-1}$ preserves $A$, we obtain
\begin{align*}
u=S_a^{-1}(S_a(u))=S_a^{-1}(p)ru.
\end{align*}
Uniqueness of the factorization gives $S_a^{-1}(p)r=1_A$.
By Lemma~\ref{lem:onesided-MR-units}, $S_a^{-1}(p)$ is a unit of $A$.
Since $S_a$ restricts to an automorphism of $A$, we also have $p\in A^\times$.
By condition~\eqref{MR1},
\begin{align*}
ua=aS_a(u)=apq.
\end{align*}
Uniqueness of the factorization gives $u\tright a=ap$ and $u\tleft a=q$.
Hence $u\tright a\in aA^\times$.

The second assertion follows similarly, using $S_u^{-1}$ in place of $S_a$.
\end{proof}

\begin{prop}\label{prop:universal-surjectivity}
For every Zappa--Sz\'ep product
$M=A\bowtie U$ of multiplicative racks, the ordered product of subgroups satisfies
\begin{align*}
\cG(M)=\cG(A)\cG(U).
\end{align*}
\end{prop}

\begin{proof}
By Proposition~\ref{prop:universal-inclusions}, we regard $\cG(A)$ and $\cG(U)$ as subgroups of $\cG(M)$, and they generate $\cG(M)$.
Let $N=\omega_M(M^\times)$.
For $x\in M$ and $h\in M^\times$, condition~\eqref{MR1} gives
\begin{align*}
\omega_M(x)^{-1}\omega_M(h)\omega_M(x) =\omega_M(h*x).
\end{align*}
Since every right translation $S_x$ is a monoid automorphism, it maps $M^\times$ bijectively onto itself.
Consequently,
\begin{align*}
\omega_M(x)^{-1}N\omega_M(x)=N
\end{align*}
for all $x\in M$.
Since $\omega_M(M)$ generates $\cG(M)$, the subgroup $N$ is normal in $\cG(M)$.

We observe that $M^\times=A^\times U^\times$.
Indeed, let $au\in M^\times$ and write $(au)^{-1}=bv$, where $a,b\in A$ and $u,v\in U$.
Comparing the components of $(au)(bv)=1$ and $(bv)(au)=1$ gives
\begin{align*}
a(u\tright b)=1_A, \quad (v\tleft a)u=1_U.
\end{align*}
By Lemma~\ref{lem:onesided-MR-units}, we obtain $a\in A^\times$ and $u\in U^\times$.
The reverse inclusion is immediate.
Consequently,
\begin{align*}
N=\omega_M(A^\times)\omega_M(U^\times).
\end{align*}
By Lemma~\ref{lem:unit-shifts}, for every $a\in A$ and $u\in U$ there exist $x\in A^\times$ and $y\in U^\times$ such that $ua=axyu$.
Hence the images of $\cG(A)$ and $\cG(U)$ commute in $\cG(M)/N$.
Since these images generate $\cG(M)/N$, every $g\in \cG(M)$ can be written as $g=hkn$ for some $h\in\cG(A)$, $k\in\cG(U)$, and $n\in N$.
Since $N$ is normal in $\cG(M)$, we have $knk^{-1}\in N$.
Thus we may write $ knk^{-1}=pq$ for some $p\in\omega_M(A^\times)$ and $q\in\omega_M(U^\times)$.
It follows that
\begin{align*}
g =hkn =h(knk^{-1})k =(hp)(qk)\in\cG(A)\cG(U).
\end{align*}
Therefore, we obtain $\cG(M)=\cG(A)\cG(U)$.
\end{proof}

It remains to establish the uniqueness of the ordered decomposition.
Lemma \ref{lem:unit-shifts} will allow us to reduce the intersection of the two subgroups to the images of units in the original factors.

\begin{thm}\label{thm:universal-exactness}
Let $M=A\bowtie U$ be a Zappa--Sz\'ep product of multiplicative racks.
The canonical homomorphisms $\cG(A) \to \cG(M)$ and $\cG(U) \to \cG(M)$ are injective, and their images form an exact factorization
\begin{align*}
\cG(M)&\cong\cG(A)\bowtie\cG(U).
\end{align*}
\end{thm}

\begin{proof}
By Propositions~\ref{prop:universal-inclusions} and~\ref{prop:universal-surjectivity}, it suffices to show $\cG(A)\cap\cG(U)=\{1\}$.
Let $g$ belong to this intersection.
Using the left- and right-fraction descriptions, write
\begin{align*}
g&=\omega_A(a)^{-1}\omega_A(b)=\omega_U(u)\omega_U(v)^{-1}
\end{align*}
for $a,b\in A$ and $u,v\in U$.
Thus $\omega_M(au)=\omega_M(bv)$.
By Lemma~\ref{lem:cancellative-quotient}, there exist $c\in A$ and $z\in U$ such that
\begin{align*}
(au)(cz) = (bv)(cz).
\end{align*}
The $U$-components give $(u \tleft c)z = (v \tleft c)z$.
By Lemma~\ref{lem:unit-shifts}, write
\begin{align*}
u \tleft c = pu, \quad v \tleft c = qv
\end{align*}
for some $p,q \in U^\times$.
Passing to $\cG(U)$ yields $g=\omega_U(p^{-1}q)$.
Put $k=p^{-1}q \in U^\times$.
Then $\omega_M(ak) = \omega_M(b)$.
By Lemma~\ref{lem:cancellative-quotient} again, there exist $d \in A$ and $w \in U$ such that
\begin{align*}
(ak)(dw) = b(dw).
\end{align*}
Condition~\eqref{MR2} and the bijectivity of $S_{dw}$ give $S_{ak} = S_b$, and hence $S_k = S_b \circ S_a^{-1}$.
Since $a,b\in A$, both $S_a$ and $S_b$ restrict to automorphisms of $A$, and hence $S_k(A)=A$.
Moreover, $k\in U^\times$ and condition~\eqref{MR2} give $S_{k^{-1}}=S_k^{-1}$, and hence $S_{k^{-1}}(d)\in A$.
By condition~\eqref{MR1}, $kd = S_k^{-1}(d)k$ holds.
Therefore
\begin{align*}
(aS_k^{-1}(d))(kw) = (bd)w.
\end{align*}
The uniqueness of the factorization $M=A \bowtie U$ gives $kw = w$ in $U$.
Applying $\omega_U$ yields $\omega_U(k)=1$, and hence $g=1$.
\end{proof}

\subsection{Reduction modulo units}\label{subsec:unit-reduction}
For a multiplicative rack $M$, define an equivalence relation by
\begin{align*}
x\sim_{\red}y \iff y=xh\ \text{for some }h\in M^\times.
\end{align*}
Put $M_{\red}\coloneq M/{\sim_{\red}}$, and let $[x]$ denote the equivalence class of $x\in M$ in this quotient.
The following lemma justifies treating this quotient as a multiplicative rack.

\begin{lemma}\label{lem:unit-reduction}
Let $M$ be a multiplicative rack.
Then the relation $\sim_{\red}$ is a congruence for both the monoid multiplication and the rack operation, and the quotient $M_{\red}$ is a multiplicative rack.
Moreover, $(M_{\red})^\times$ is trivial.
\end{lemma}

\begin{proof}
Since $M^\times$ is a group, $\sim_{\red}$ is an equivalence relation.
We first show that it is compatible with multiplication.
As observed in the proof of Proposition~\ref{prop:universal-surjectivity}, every right translation preserves $M^\times$.
Thus, for $h,k\in M^\times$ and $x,y\in M$, condition~\eqref{MR1} gives
\begin{align*}
(xh)(yk)=xy(h*y)k\sim_{\red}xy.
\end{align*}
Hence $\sim_{\red}$ is a monoid congruence.
It is also compatible with the rack operation in the first variable, since condition~\eqref{MR3} gives
\begin{align*}
(xh)*z=(x*z)(h*z)\sim_{\red}x*z.
\end{align*}
For the second variable, condition~\eqref{MR1} implies
\begin{align*}
S_k(t)=k^{-1}tk=t\bigl((k^{-1}*t)k\bigr) \sim_{\red} t \quad(t\in M),
\end{align*}
where $(k^{-1}*t)k$ is a unit.
Thus, by condition~\eqref{MR2},
\begin{align*}
x*(zk)=S_k(S_z(x))\sim_{\red}x*z.
\end{align*}
Since $S_z$ and $S_z^{-1}$ preserve $M^\times$, the induced right translations on the quotient are bijective.
All other axioms of the multiplicative rack descend from $M$.

Finally, if $[x]$ is a unit in $M_{\red}$, choose $y\in M$ with $[xy]=[1]$.
Then $xy=h\in M^\times$, so $x(yh^{-1})=1$.
By Lemma~\ref{lem:onesided-MR-units}, $x\in M^\times$, and therefore $[x]=[1]$.
\end{proof}

\begin{thm}\label{thm:unit-reduction-ZS}
Let $M=A\bowtie U$ be a Zappa--Sz\'ep product of multiplicative racks.
Then the groups of units form an exact factorization of groups
\begin{align*}
M^\times=A^\times\bowtie U^\times,
\end{align*}
and there is a canonical isomorphism of multiplicative racks
\begin{align*}
M_{\red}\cong A_{\red}\times U_{\red}.
\end{align*}
\end{thm}

\begin{proof}
As shown in the proof of Proposition~\ref{prop:universal-surjectivity}, we have
\begin{align*}
M^\times=A^\times U^\times.
\end{align*}
The uniqueness of the factorization $M=A\bowtie U$ makes this an exact factorization of groups.

It remains to establish the asserted isomorphism of multiplicative racks.
Consider the map
\begin{align*}
\varphi \colon A_{\red}\times U_{\red} \to M_{\red};\ ([a],[u])\mapsto[au].
\end{align*}
The map is well-defined: if $a'=ap$ and $u'=uq$ for some $p\in A^\times$ and $q\in U^\times$, then condition~\eqref{MR1} gives
\begin{align*}
a'u'=apuq=au(p*u)q\sim_{\red}au,
\end{align*}
since $p*u\in M^\times$.
By Lemma~\ref{lem:unit-shifts}, for every $a\in A$ and $u\in U$, there exist $p\in A^\times$ and $q\in U^\times$ such that $ua=apqu$.
Applying condition~\eqref{MR1}, we obtain
\begin{align*}
ua=au((pq)*u)\sim_{\red}au.
\end{align*}
Thus the images of $A$ and $U$ commute in $M_{\red}$, which shows that $\varphi$ is a monoid homomorphism.
It is surjective because $M=AU$.

To prove injectivity, suppose $[au]=[bv]$ in $M_{\red}$.
Then $bv=auh$ for some $h\in M^\times$.
Write $h=pq$ with $p\in A^\times$ and $q\in U^\times$.
The defining Zappa--Sz\'ep relation gives
\begin{align*}
bv=a(u\tright p)(u\tleft p)q.
\end{align*}
By Lemma~\ref{lem:unit-shifts}, $u\tright p\in pA^\times=A^\times$ and $u\tleft p=tu$ for some $t\in U^\times$.
Uniqueness of the factorization gives
\begin{align*}
b=a(u\tright p),\quad v=tuq=u(t*u)q,
\end{align*}
where $t*u\in U^\times$.
Hence $[a]=[b]$ in $A_{\red}$ and $[u]=[v]$ in $U_{\red}$.
So $\varphi$ is a monoid isomorphism.

By Lemma~\ref{lem:unit-reduction}, all three quotients are multiplicative racks, and the quotient maps are homomorphisms of multiplicative racks.
The monoid isomorphism $\varphi$ identifies the images of $A$ and $U$ in $M_{\red}$ with $A_{\red}$ and $U_{\red}$, respectively; these images are multiplicative subracks.
It remains to verify that the rack operation is componentwise.
For $a\in A$ and $u\in U$, the proof of Lemma~\ref{lem:unit-shifts} gives
\begin{align*}
S_a(u)=p(u\tleft a),\quad S_u^{-1}(a)=(u\tright a)q
\end{align*}
for some $p\in A^\times$ and $q\in U^\times$.
We may also write $u\tleft a=tu$ and $u\tright a=ar$ for some $t\in U^\times$ and $r\in A^\times$.
Thus condition~\eqref{MR1} gives
\begin{align*}
S_a(u)=ptu=u((pt)*u),\quad S_u^{-1}(a)=arq.
\end{align*}
Since $(pt)*u$ and $rq$ are units of $M$, we obtain
\begin{align*}
[S_a(u)]=[u],\quad [S_u^{-1}(a)]=[a].
\end{align*}
Because the induced right translations are bijective, the second equality also implies $[S_u(a)]=[a]$.
Thus the induced right translations of either factor fix the image of the other factor pointwise.
Using conditions~\eqref{MR2} and \eqref{MR3}, we conclude that
\begin{align*}
[(au)*(bv)]&=[S_v(S_b(au))]\\
&=[((a*b)(u*b))*v]\\
&=[((a*b)*v)((u*b)*v)]\\
&=[a*b][u*v]
\end{align*}
for $a,b\in A$ and $u,v\in U$.
Consequently $\varphi$ preserves the componentwise rack operation and is an isomorphism of multiplicative racks.
\end{proof}

\begin{cor}\label{cor:trivial-units-direct}
Let $M=A\bowtie U$ be a Zappa--Sz\'ep product of multiplicative racks.
If $A^\times=\{1_A\}$ $($resp.\ $U^\times=\{1_U\})$, then $u\tright a=a$ $($resp.\ $u\tleft a=u)$ for all $a\in A$ and $u\in U$.
In particular, if both groups of units are trivial, then $M$ is canonically isomorphic to the direct product $A\times U$ as a multiplicative rack.
\end{cor}

\begin{proof}
The first assertion follows immediately from Lemma~\ref{lem:unit-shifts}.
If both groups of units are trivial, then $A_{\red}=A$ and $U_{\red}=U$, so the second assertion follows from Theorem~\ref{thm:unit-reduction-ZS}.
\end{proof}

\subsection{Cancellative quotients and Ore localization}\label{subsec:localization-applications}
We now derive two consequences of Theorem~\ref{thm:universal-exactness} concerning cancellative quotients and Ore localization.

\begin{cor}\label{cor:canonical-factorization}
Let $M=A\bowtie U$ be a Zappa--Sz\'ep product of multiplicative racks.
For all $a,b\in A$ and $u,v\in U$, the following equivalence holds$:$
\begin{align*}
\omega_M(au)=\omega_M(bv) \iff \omega_A(a)=\omega_A(b)\ \text{and}\ \omega_U(u)=\omega_U(v).
\end{align*}
Consequently, there is a canonical exact factorization of multiplicative racks
\begin{align*}
M_{\can}&\cong A_{\can}\bowtie U_{\can},
\end{align*}
whose actions, given by $L_{[u]}([a])=[u\tright a]$ and $R_{[a]}([u])=[u\tleft a]$, are bijective.
\end{cor}
\begin{proof}
By Theorem~\ref{thm:universal-exactness}, the universal groups of the factors form an exact factorization of $\cG(M)$.
This gives the equivalence above by uniqueness of the ordered decomposition.
Consequently, the multiplication map
\begin{align*}
A_{\can}\times U_{\can} \to M_{\can};\ ([a],[u]) \mapsto [au]
\end{align*}
is bijective.
The factors $A_{\can}$ and $U_{\can}$ embed into $M_{\can}$ as multiplicative subracks, and thus we obtain $M_{\can}\cong A_{\can}\bowtie U_{\can}$ as multiplicative racks.
The induced Zappa--Sz\'ep action is given by
\begin{align*}
[u]\tright[a]=[u\tright a],\quad [u]\tleft[a]=[u\tleft a].
\end{align*}
These formulas follow from the defining relation $ua=(u\tright a)(u\tleft a)$ and uniqueness of the factorization in $M_{\can}$.
The corresponding coordinate maps are bijective by Proposition~\ref{prop:MR ZSaction bijective}.
\end{proof}

\begin{cor}\label{cor:localization}
Let $M=A\bowtie U$ be a Zappa--Sz\'ep product of multiplicative racks.
If $M$ is left-cancellative (equivalently, cancellative), then its groups of right fractions admit a canonical exact factorization
\begin{align*}
\what M&\cong\what A\bowtie\what U.
\end{align*}
\end{cor}
\begin{proof}
The submonoids $A,U$ are cancellative, and Lemma~\ref{lemma:MRhasCM} shows that all three monoids are Ore.
Their universal groups are therefore the corresponding groups of right fractions, so the result follows from Theorem~\ref{thm:universal-exactness}.
\end{proof}

The analogous statement fails for general Ore monoids.
For instance, consider the positive Baumslag--Solitar monoid
\begin{align*}
M=\langle a,b\mid ba=ab^2\rangle_{\Mon}.
\end{align*}
Every element of $M$ has a unique expression $a^ib^j$ with $i,j\ge0$;
see for example \cite[Section 3.1]{BrownloweRamaggeRobertsonWhittaker2014}.
Thus $M=A\bowtie U$, where $A=\langle a\rangle$ and $U=\langle b\rangle$.
The multiplication rule
\begin{align*}
(a^ib^j)(a^kb^\ell) =a^{i+k}b^{2^kj+\ell}
\end{align*}
shows that $M$ is cancellative and has common right multiples.
Thus $M$ is an Ore monoid.
However, the ordered product $\cG(A)\cG(U)$ does not exhaust $\cG(M)$.
Indeed, suppose that $aba^{-1}=a^mb^n$ for some $m,n\in\bZ$.
The homomorphism $\cG(M)\to\bZ$ sending $a$ to $1$ and $b$ to $0$ gives $m=0$.
Thus $aba^{-1}=b^n$, and hence $b=ab^2a^{-1}=b^{2n}$, contradicting the fact that $b$ has infinite order, since $M$ embeds into $\cG(M)$.
Therefore, we have $\cG(A)\cG(U)\subsetneq\cG(M)$.

\section{Unital ideals and generalized Bianchi-type factorizations}\label{sec:non-cancellative}
In this section, we study multiplicative racks with unital ideals and their Zappa--Sz\'ep factorizations.
We first generalize Bianchi's disjoint-union construction~\cite{Bianchi2025} by allowing augmentations into arbitrary multiplicative racks rather than groups.
We then establish a reconstruction theorem for multiplicative racks with suitable unital ideals and characterize their exact factorizations.

\subsection{Partially multiplicative racks and generalized Bianchi-type construction}\label{subsec:PMQ}
We first recall Bianchi's notion of a partially multiplicative rack~\cite[Appendix]{Bianchi2025}.
We then introduce the data needed to generalize his disjoint-union construction.

\begin{defi}
Let $Q$ be a unital rack with unit $1$.
Suppose that $Q$ is equipped with a partial multiplication
\begin{align*}
D \to Q;\ (a,b) \mapsto ab,
\end{align*}
where $D \subseteq Q \times Q$.

Following Bianchi's partial-monoid convention~\cite[Definition~2.3 and Appendix]{Bianchi2025}, we say that $Q$ is a \textbf{partially multiplicative rack} if the following conditions hold:
\begin{enumerate}
\conditem{PM1}\label{PM1} the products $1a$ and $a1$ are always defined and equal to $a$,
\conditem{PM2}\label{PM2} the products $ab$ and $(ab)c$ are both defined if and only if $bc$ and $a(bc)$ are both defined, and in this case $(ab)c=a(bc)$,
\conditem{PM3}\label{PM3} the product $ab$ is defined if and only if $b(a*b)$ is defined, and in this case $ab=b(a*b)$,
\conditem{PM4}\label{PM4} if $bc$ is defined, then $a*(bc)=(a*b)*c$,
\conditem{PM5}\label{PM5} the product $ab$ is defined if and only if $(a*c)(b*c)$ is defined, and in this case $(ab)*c=(a*c)(b*c)$.
\end{enumerate}
A partially multiplicative rack is called a \textbf{partially multiplicative quandle} if its underlying rack is a quandle.
A partially multiplicative rack (resp. partially multiplicative quandle) is \textbf{complete} if $D = Q \times Q$.
\end{defi}

We introduce the following generalization of Bianchi's PMR-group pair.

\begin{defi}[cf. {\cite[Definition 2.15]{Bianchi2025}}]
A \textbf{PMR--MR pair} is a quadruple $(Q,I,\varepsilon,\rho)$ consisting of a partially multiplicative rack $Q$, a multiplicative rack $I$, a homomorphism of partially multiplicative racks $\varepsilon\colon Q\to I$, and a right monoid action $\rho$ of $I$ on $Q$ by automorphisms of partially multiplicative racks.
We call $\varepsilon$ the \textbf{augmentation} and write $x^i\coloneq\rho_i(x)$.
These data are required to satisfy the following conditions:
\begin{enumerate}
\conditem{GP1}\label{GP1} $x^{\varepsilon(y)} = x*y \quad (x,y \in Q)$,
\conditem{GP2}\label{GP2} $\varepsilon(x^i) = \varepsilon(x)*i \quad (x \in Q, i \in I)$.
\end{enumerate}
When $Q$ is a partially multiplicative quandle, we call the structure a \textbf{PMQ--MR pair}.
\end{defi}

\begin{prop}[cf. {\cite[Definition 2.21]{Bianchi2025}}] \label{prop:PMR-completion}
Let $(Q, I, \varepsilon, \rho)$ be a PMR--MR pair.
Consider the disjoint union $\cB(Q,I) = Q \sqcup \overline{I}$, where $\overline{I} = \left\{ \overline{i} \mid i \in I \right\}$ is a copy of $I$.
Define a monoid operation on $\cB(Q,I)$ by retaining the partial multiplication on $Q$ when defined and the multiplication on $\overline I$ induced from $I$, and by
\begin{align*}
&xy \coloneq \overline{\varepsilon(x) \varepsilon(y)}\ \text{if the product of $x,y \in Q$ is not defined in $Q$},\\
&x\overline{i}\coloneq \overline{\varepsilon(x) i},\\
&\overline{i}x \coloneq \overline{i\varepsilon(x)},
\end{align*}
for $x,y \in Q$ and $i \in I$.
Define a rack operation on $\cB(Q,I)$ by retaining the rack operation on $Q$ and $I$, and by
\begin{align*}
&x*\overline{i} \coloneq x^i,\\
&\overline{i}*x \coloneq \overline{i*\varepsilon(x)},
\end{align*}
for $x \in Q$ and $i \in I$.
Then $\cB(Q,I)$ is a multiplicative rack.
Moreover, $\cB(Q,I)$ is a multiplicative quandle whenever $Q$ is a partially multiplicative quandle and $I$ is a multiplicative quandle.
\end{prop}

\begin{proof}
We first verify that the multiplication on $\cB(Q,I)$ is associative.
For $x,y,z\in Q$, partial associativity condition~\eqref{PM2} shows that $(xy)z$ and $x(yz)$ both lie in $Q$ precisely when both successive products are defined in $Q$; in that case they agree by condition~\eqref{PM2}.
In every other case both results lie in $\overline I$ and agree because $\varepsilon$ preserves the defined products and multiplication in $I$ is associative.
If at least one of the three factors belongs to $\overline I$, associativity follows directly from the mixed multiplication rules and the multiplicativity of $\varepsilon$ on defined products.
The element $1_Q$ is a two-sided identity by condition~\eqref{PM1} and $\varepsilon(1_Q)=1_I$.
Thus $\cB(Q,I)$ is a monoid.
Define a map $\mu \colon \cB(Q,I) \to I$ by
\begin{align*}
\mu(x)&\coloneq \varepsilon(x) \quad (x \in Q),\\
\mu(\overline{i}) &\coloneq i \quad (i \in I).
\end{align*}
The multiplication rules show that $\mu$ is a monoid homomorphism.
For each $i \in I$, we define a map $T_i \colon \cB(Q,I) \to \cB(Q,I)$ by
\begin{align*}
T_i(x) &\coloneq x^i = \rho_i(x) \quad (x \in Q),\\
T_i(\overline{j}) &\coloneq \overline{j*i} \quad (j \in I).
\end{align*}
Since $\rho_i$ is an automorphism of $Q$ as a partially multiplicative rack and the right translation on $I$ is an automorphism of the multiplicative rack, the map $T_i$ is bijective.
By the definition of $T_i$, we have $T_{ij} = T_j \circ T_i$ for any $i,j \in I$.
By condition~\eqref{GP1} and the definition of the rack operation, we have $S_q=T_{\mu(q)}$ for every $q\in\cB(Q,I)$.
We verify the remaining conditions of Remark~\ref{remark:MR}.
First, we show that $pq=qS_q(p)$ for all $p,q\in\cB(Q,I)$.
If $x,y\in Q$ and $xy$ is defined in $Q$, this follows from condition~\eqref{PM3}.
If $xy$ is not defined in $Q$, then $y(x*y)$ is also undefined in $Q$, and
\begin{align*}
y(x*y) = \overline{\varepsilon(y) \varepsilon(x*y)} = \overline{\varepsilon(y) (\varepsilon(x)*\varepsilon(y))} = \overline{\varepsilon(x) \varepsilon(y)} = xy.
\end{align*}
For $x \in Q$ and $i \in I$, we have
\begin{align*}
\overline{i} (x*\overline{i}) &= \overline{i\varepsilon(x^i)} = \overline{i(\varepsilon(x)*i)} = \overline{\varepsilon(x)i} = x\overline{i},\\
x(\overline{i}*x) &= \overline{\varepsilon(x) (i*\varepsilon(x))} = \overline{i\varepsilon(x)} = \overline{i}x.
\end{align*}
Next, since $\mu$ is a monoid homomorphism, we obtain
\begin{align*}
S_{qr} = T_{\mu(qr)} = T_{\mu(q)\mu(r)} = T_{\mu(r)} \circ T_{\mu(q)} = S_r \circ S_q.
\end{align*}
Finally, we show that each $T_i$ is a monoid automorphism.
Since $\rho_i$ is an automorphism of the partially multiplicative rack $Q$, the map $T_i$ preserves products defined in $Q$.
For $x,y\in Q$ whose product is undefined in $Q$, conditions~\eqref{GP2} and \eqref{MR3} in $I$ give
\begin{align*}
T_i(xy)=\overline{(\varepsilon(x)\varepsilon(y))*i}=\overline{(\varepsilon(x)*i)(\varepsilon(y)*i)}=T_i(x)T_i(y).
\end{align*}
The mixed products and products in $\overline I$ are preserved by the same compatibility identities.
Since $T_i$ is bijective and preserves multiplication, it is a monoid automorphism.
Thus Remark~\ref{remark:MR} shows that $\cB(Q,I)$ is a multiplicative rack.
Moreover, $\cB(Q,I)$ is a multiplicative quandle if and only if $Q$ is a partially multiplicative quandle and $I$ is a multiplicative quandle.
\end{proof}

\begin{remark}
When $I$ is a group $G$, regarded as a multiplicative quandle under conjugation, a PMR--MR pair reduces to a PMR-group pair in Bianchi's setting.
Moreover, the generalized Bianchi-type construction recovers Bianchi's disjoint-union construction.
Bianchi uses the notation $Q\bowtie G$ for this disjoint-union construction.
We write $\cB(Q,G)$ to distinguish it from a Zappa--Sz\'ep product of monoids.
Note that the underlying monoid of the Bianchi-type construction $Q\sqcup \overline{G}$ is always non-cancellative.
Indeed, the distinct elements $1_Q$ and $\overline{1_G}$ satisfy $1_Q \overline{1_G} = \overline{1_G} = \overline{1_G}\,\overline{1_G}$.
\end{remark}

\subsection{Unital ideals and generalized Bianchi-type presentations}\label{subsec:ideal}
We study unital ideals of multiplicative racks and their role in generalized Bianchi-type presentations.
After establishing basic properties of unital and minimal ideals, we prove a reconstruction theorem and discuss its consequences.

Recall that a nonempty subset $I$ of a monoid $M$ is a (two-sided) \textbf{ideal} if
\begin{align*}
MI \subseteq I \quad\text{and}\quad IM \subseteq I.
\end{align*}
An ideal is said to be \textbf{unital} if it has the identity element.
In particular, it is a monoid.
If it forms a group under the multiplication inherited from $M$, then it is called a \textbf{group ideal}.
For example, in Bianchi's construction, $\overline{G}$ is a group ideal in $\cB(M,G)$.

An ideal is called \textbf{minimal} if it contains no proper ideal of $M$.
A minimal ideal of a monoid is unique whenever it exists, since any two ideals have nonempty intersection.
Moreover, a group ideal is automatically the unique minimal ideal.
Indeed, let $I$ be a group ideal with identity $e$, and $J$ be an ideal.
For any $x\in J$, we have $xe\in I\cap J$.
Since $I$ is a group, we obtain $I = I(xe) \subseteq J$.

\begin{lemma}\label{lemma:unital ideal}
Let $M$ be a monoid, and $I$ be a unital ideal of $M$.
Then the following hold$:$
\begin{enumerate}
\item If $e$ is the identity of $I$, then $e$ is a central idempotent in $M$.
\item The map
\begin{align*}
\mu \colon M \to I; x \mapsto xe
\end{align*}
is a surjective monoid homomorphism such that $\mu|_I = \id_I$.
\item The universal group $\cG(M)$ is isomorphic to $\cG(I)$.
\end{enumerate}
\end{lemma}

\begin{proof}
(1) Since $I$ is an ideal, both $ex$ and $xe$ belong to $I$ for every $x\in M$.
Hence
\begin{align*}
ex=(ex)e=e(xe)=xe.
\end{align*}
Thus $e$ is a central idempotent of $M$.\\
(2) We see that $\mu$ is a monoid homomorphism by (1).
Since $ie=i$ for every $i\in I$, we have $\mu|_I=\id_I$.\\
(3) Let $H$ be a group and let $f\colon M\to H$ be a monoid homomorphism.
Since $e^2=e$, we have $f(e)=1_H$.
Consequently,
\begin{align*}
f(x)=f(x)f(e)=f(xe)=f(\mu(x))
\end{align*}
for every $x\in M$.
Thus $f$ factors uniquely through $\mu$.
By the universal property of universal groups, $\mu$ induces a canonical isomorphism $\mathcal G(M)\cong\mathcal G(I)$.
\end{proof}

We call the map $\mu \colon M \to I$ defined in Lemma~\ref{lemma:unital ideal} the \textbf{retraction}.

\begin{prop} \label{prop:minimal ideal}
Every minimal ideal of a multiplicative rack is a group.
\end{prop}

\begin{proof}
Suppose that $I$ is a minimal ideal of a multiplicative rack $M$.
For $a \in I$, Lemma~\ref{lemma:MRidempotent} gives $Ma=aM$.
Thus
\begin{align*}
M(aI)=(Ma)I\subseteq aI \quad \text{and} \quad (aI)M=a(IM)\subseteq aI,
\end{align*}
and similarly
\begin{align*}
M(Ia)\subseteq Ia \quad \text{and} \quad (Ia)M=I(aM)=I(Ma)\subseteq Ia.
\end{align*}
Hence $aI$ and $Ia$ are ideals contained in $I$.
By minimality, we have $aI = Ia = I$.
There exists $e_R, e_L \in I$ such that $ae_R = a$ and $e_La = a$.
For any $x\in I$, we can write $x = ba = ac$ for some $b,c \in I$.
Hence
\begin{align*}
xe_R=bae_R=ba=x, \quad e_Lx=e_Lac=ac=x.
\end{align*}
Thus $e_R$ and $e_L$ are a right and a left identity of $I$, respectively.
Since $e_L=e_Le_R=e_R$, their common value $e$ is the identity of $I$.
For each $a \in I$, there exist $b,c \in I$ such that $ab = e = ca$, and then $b = (ca)b = c(ab) = c$.
Thus $b=c$ is the inverse of $a$ in $I$.
Therefore $I$ is a group.
\end{proof}

\begin{cor}\label{cor:MR minimal ideal}
Let $M$ be a multiplicative rack having a unital ideal $I$.
Then $S_x = S_{\mu(x)}$ holds for any $x \in M$, where $\mu \colon M \to I$ is the retraction.
If $I$ is minimal, there is a canonical isomorphism $\cG(M) \cong I$.
\end{cor}

\begin{proof}
Let $e$ be the identity of $I$.
Since $e^2=e$, Lemma~\ref{lemma:MRidempotent} gives $S_e=\id_M$.
Hence, by condition~\eqref{MR2},
\begin{align*}
S_{\mu(x)} =S_{xe} =S_e\circ S_x =S_x.
\end{align*}
The final assertion follows from Lemma~\ref{lemma:unital ideal}, since $I$ is a group by Proposition~\ref{prop:minimal ideal}.
\end{proof}

If the minimal ideal $I$ is proper, then its identity $e$ satisfies $e\neq 1_M$ and $1_Me=e^2=e$.
Hence $M$ is non-cancellative.
A proper minimal ideal determines the universal group of a multiplicative rack.
The following example illustrates this phenomenon.

\begin{ex}[Lebed's construction]
The following multiplicative quandle was considered by Lebed~\cite{Lebed2014}.
For a commutative monoid $X$ with zero and a positive integer $n$, let
\begin{align*}
Q_{X,n}&=\{(\bm{x},\sigma)\in X^n\times S_n \mid x_i=0\text{ whenever }\sigma(i)\ne i\}.
\end{align*}
With coordinatewise multiplication $\bm{x}\bm{y}$ and the usual right action of $S_n$ on $X^n$, define
\begin{align*}
(\bm{x},\sigma)(\bm{y},\tau)& \coloneq (\bm{x}\bm{y},\sigma\tau),\\
(\bm{x},\sigma)*(\bm{y},\tau)& \coloneq (\bm{x}\tau,\tau^{-1}\sigma\tau).
\end{align*}
For nontrivial $X$, this multiplicative quandle is non-cancellative.
Its minimal group ideal is
\begin{align*}
\{(\bm{0},\sigma) \mid \sigma\in S_n\}\cong S_n.
\end{align*}
Indeed, multiplication by $(\bm{0},\sigma^{-1}\tau)$ sends $(\bm{x},\sigma)$ to $(\bm{0},\tau)$.
Thus Corollary~\ref{cor:MR minimal ideal} gives $\cG(Q_{X,n})\cong S_n$.
\end{ex}

We now establish a converse to the generalized Bianchi-type construction.
The key observation is that a unital ideal determines a canonical retraction, while invariance under right translations allows the complementary subset to inherit a partially multiplicative rack structure.

\begin{thm}\label{thm:MR-Bianchi}
Let $M$ be a multiplicative rack with a proper unital ideal $I$.
Suppose that $S_x(I) = I$ for every $x \in M$.
Then $I$ is a multiplicative rack with the induced operations.

Put $Q=M\setminus I$, and equip $Q$ with the restricted rack operation and the partial multiplication inherited from $M$, where the product of $x,y\in Q$ is defined precisely when $xy\in Q$.
Then $Q$ is a partially multiplicative rack, and
\begin{align*}
(Q, I, \mu|_Q, \rho), \quad \rho_i(x) \coloneq x*i\quad (i \in I, x \in Q)
\end{align*}
is a PMR--MR pair, where $\mu \colon M \to I$ is the retraction.

Moreover, there is a canonical isomorphism of multiplicative racks $M \cong \cB(Q,I)$.
In particular, $M$ is a multiplicative quandle if and only if $Q$ is a partially multiplicative quandle and $I$ is a multiplicative quandle.
\end{thm}

\begin{proof}
First, since $I$ is a unital ideal, it is a monoid with identity $e$.
By Lemma~\ref{lemma:unital ideal}, $e$ is idempotent, and hence $S_e=\id_M$ by Lemma~\ref{lemma:MRidempotent}.
Moreover, since $S_x(I)=I$ and $S_x$ is a monoid automorphism for every $x\in M$, we have $S_x(e)=e$.
Thus $e$ is a rack unit of $I$, and each $S_i$ restricts to a monoid automorphism of $I$ for $i\in I$.
Therefore, $I$ is a multiplicative rack.

Next, we show that $Q = M \setminus I$ is a partially multiplicative rack.
Since $I$ is a proper ideal, we have $1_M \in Q$.
By assumption, $S_x(Q)=Q$ for every $x \in M$, so $Q$ inherits a unital rack structure from $M$.
The partial multiplication on $Q$ satisfies the unit and partial associativity axioms.
Indeed, since $I$ is an ideal, for any $x,y,z \in Q$, the product $(xy)z$ is defined in $Q$ if and only if $xyz \in Q$.
The same holds for $x(yz)$.
Thus partial associativity follows from associativity in $M$.
The identities required by conditions~\eqref{PM3}--\eqref{PM5}, whenever the corresponding products are defined, follow from conditions~\eqref{MR1}--\eqref{MR3}.
To verify the domain conditions, observe that for all $x,y,z\in Q$, we have
\begin{align*}
xy \in Q \iff y(x*y) \in Q
\end{align*}
by condition~\eqref{MR1}, and
\begin{align*}
xy \in Q \iff (xy)*z \in Q \iff (x*z)(y*z) \in Q
\end{align*}
by condition~\eqref{MR3} and $S_z(Q)=Q$.
Thus $Q$ is a partially multiplicative rack.

We verify that $(Q,I,\mu|_Q,\rho)$ is a PMR--MR pair.
First, the map $\mu$ is a homomorphism of multiplicative racks.
Indeed, by Lemma~\ref{lemma:unital ideal}, it is a monoid homomorphism.
Moreover, since $S_y(e) = e$ for every $y \in M$, we have
\begin{align*}
\mu(x*y) = (x*y)e = (x*y)(e*y) = (xe)*y = (xe)*(ye) = \mu(x)*\mu(y).
\end{align*}
Therefore $\mu|_Q \colon Q \to I$ is a homomorphism of partially multiplicative racks.
For each $i\in I$, the restriction $S_i|_Q$ is an automorphism of $Q$ as a partially multiplicative rack.
Since $S_e=\id_M$, condition~\eqref{MR2} shows that $\rho$ defines a right action of $I$ on $Q$ by such automorphisms.
By Corollary~\ref{cor:MR minimal ideal}, $S_{\mu(x)} = S_x$ holds, and thus
\begin{align*}
x^{\mu(y)} = x*\mu(y) = x*y
\end{align*}
for all $x,y\in Q$, where $x^i \coloneq \rho_i(x)$.
Furthermore, for $x \in Q$ and $i\in I$,
\begin{align*}
\mu(x^i) = \mu(x*i) = \mu(x)*\mu(i) = \mu(x)*i.
\end{align*}
Thus the augmentation is equivariant.
Therefore $(Q,I,\mu|_Q,\rho)$ is a PMR--MR pair.

Finally, we show that the map $\varphi \colon \cB(Q,I) \to M$ defined by
\begin{align*}
\varphi(x) \coloneq x \quad (x \in Q), \quad \varphi(\overline{i}) \coloneq i \quad (i \in I)
\end{align*}
preserves the multiplication and rack operation.
Note that $\cB(Q,I) = Q \sqcup \overline{I}$ can be identified with $M$ as sets, and thus $\varphi$ is bijective.
For $x,y \in Q$ whose product is defined in $Q$, this is immediate.
If their product is undefined in $Q$, then $xy \in I$, and hence
\begin{align*}
\varphi(xy) = \varphi(\overline{\mu(x)\mu(y)}) = \varphi(\overline{\mu(xy)}) = \mu(xy) = xy = \varphi(x)\varphi(y).
\end{align*}
For $x\in Q$ and $i\in I$, we have
\begin{align*}
&\varphi(x\overline{i}) = \varphi(\overline{\mu(x)i}) = \mu(x)i = (xe)i = xi = \varphi(x)\varphi(\overline{i}),\\
&\varphi(\overline{i}x) = \varphi(\overline{i\mu(x)}) = i\mu(x) = i(xe) = i(ex) = ix = \varphi(\overline{i})\varphi(x).
\end{align*}
Thus $\varphi$ is a monoid isomorphism.
It also preserves the rack operation.
Indeed, for $x\in Q$ and $i\in I$, we have
\begin{align*}
&\varphi(x*\overline{i}) = \varphi(\rho_i(x)) = \rho_i(x) = x*i = \varphi(x)*\varphi(\overline{i}),\\
&\varphi(\overline{i}*x) = \varphi(\overline{i*\mu(x)}) = i*\mu(x) = i*x = \varphi(\overline{i})*\varphi(x).
\end{align*}
The operations on $Q$ and $I$ are preserved by construction.
Thus $\varphi$ is an isomorphism of multiplicative racks.

Since $M=Q\sqcup I$ as sets, $M$ is a multiplicative quandle if and only if $Q$ is a partially multiplicative quandle and $I$ is a multiplicative quandle.
\end{proof}

Since every minimal ideal of a multiplicative rack is a group and is preserved by all automorphisms, the preceding theorem applies in particular to multiplicative racks with proper minimal ideals.

\begin{cor}\label{cor:MR-Bianchi}
Let $M$ be a multiplicative rack with a proper minimal ideal $I$.
Put $Q=M\setminus I$, and equip $Q$ with the restricted rack operation and the partial multiplication inherited from $M$, where the product of $x,y\in Q$ is defined precisely when $xy\in Q$.
Then $Q$ is a partially multiplicative rack, and
\begin{align*}
(Q, I, \mu|_Q, \rho), \quad \rho_i(x) \coloneq x*i\quad (i \in I, x \in Q)
\end{align*}
is a PMR--MR pair (in fact, a PMR-group pair in the sense of Bianchi).
Moreover, there is a canonical isomorphism of multiplicative racks $M \cong \cB(Q,I)$.
\end{cor}

\begin{proof}
By Proposition~\ref{prop:minimal ideal}, $I$ is a group.
In particular, $I$ is a unital ideal of $M$.
By Theorem~\ref{thm:MR-Bianchi}, $Q$ is a partially multiplicative rack, $(Q, I, \mu|_Q, \rho)$ is a PMR--MR pair, and $M \cong \cB(Q,I)$.
\end{proof}

\begin{remark}
Every finite monoid which is not group has a proper minimal ideal.
In particular, any finite non-group multiplicative rack admits a Bianchi-type presentation by Corollary~\ref{cor:MR-Bianchi}.
\end{remark}

The reconstruction theorem retains both the rack operation and the monoid multiplication.
This distinction is important because the underlying rack alone need not determine the multiplicative rack structure.

\begin{remark}\label{rem:coset-presentation}
Let $M$ be a multiplicative rack with minimal group ideal $I$, and let $e$ be the identity of $I$.
By Corollary~\ref{cor:MR minimal ideal}, the map $\mu:M\to I$ defined by $\mu(x)=xe$ satisfies $S_x=S_{\mu(x)}$.
Thus the right action of $I$ on $M$ given by $x^g=x*g$ determines the rack operation through $x*y=x^{\mu(y)}$.

Choose representatives $a_j$ of the $I$-orbits in $M$, and put $H_j  \coloneq \Stab_I(a_j)$ and $z_j \coloneq \mu(a_j)$.
By equivariance of $\mu$, we have $z_j\in C_I(H_j)$, where $C_I(H_j)$ denotes the centralizer of $H_j$ in $I$.
Identifying the orbit of $a_j$ with $H_j\backslash I$ via $H_jg\mapsto a_j^g$, we obtain a coset presentation of the underlying rack:
\begin{align*}
M \cong\bigsqcup_j H_j\backslash I, \quad (H_i g)*(H_j h) =H_i gh^{-1}z_jh.
\end{align*}
Moreover, $M$ is a multiplicative quandle if and only if $z_j\in H_j$ for every $j$, equivalently, $z_j\in Z(H_j)$.

This is an instance of the classical coset presentation of racks and quandles~\cite[Theorem 7.2]{Joyce1982}, with the minimal group ideal $I$ as the acting group.
However, the coset presentation determines only the rack operation, not the monoid multiplication.
For example, on $\{1,a,0\}$ with the trivial quandle operation, monoid identity $1$, and zero element $0$, the two multiplication rules $a^2=a$ and $a^2=0$ define nonisomorphic multiplicative quandles with the same underlying rack and minimal group ideal $\{0\}$.
In contrast, the Bianchi-type presentation in Corollary~\ref{cor:MR-Bianchi} retains the partial multiplication outside the minimal ideal and hence recovers the full multiplicative rack structure.
\end{remark}

\subsection{Exact factorizations and unital ideals}\label{subsec:ideal-factorizations}
We investigate how exact factorizations of monoids induce exact factorizations of their unital ideals.
We then apply these results to multiplicative racks and examine the relationship between idempotents and unital ideals meeting the factors.

\begin{prop}\label{prop:unital-ideal-factors}
Let $M=A\bowtie U$ be a Zappa--Sz\'ep product of monoids, and suppose that $M$ has a unital ideal $I$ with identity $e$.
Then the following hold$:$
\begin{enumerate}
\item If $A\cap I\neq\varnothing$, then $e \in A$ and there is an exact factorization $I = (A \cap I) \bowtie eU$, where $eU \cong U$ as monoids.
\item If $U\cap I\neq\varnothing$, then $e \in U$ and there is an exact factorization $I = eA \bowtie (U \cap I)$, where $eA \cong A$ as monoids.
\end{enumerate}
In particular, if both $A \cap I$ and $U \cap I$ are nonempty, then $I=M$.
\end{prop}

\begin{proof}
(1)\ Take $h\in A\cap I$, and write $e=au$ with $a\in A$ and $u\in U$.
Since $e$ is the identity of $I$, we have $h=he=(ha)u$.
Uniqueness of the factorization gives $u=1_U$, and hence $e=a\in A$.
For any $k=bv\in I$ with $b\in A$ and $v\in U$, we have $k=ek=(eb)v$ and $eb\in A\cap I$.
Thus $I=(A\cap I)U$, and the multiplication map
\begin{align*}
(A\cap I)\times U\to I;\ (b,v)\mapsto bv
\end{align*}
is bijective by uniqueness of the factorization of $M$.

By Lemma~\ref{lemma:unital ideal}, $e$ is a central idempotent, so the map
\begin{align*}
U\to eU;\ u\mapsto eu
\end{align*}
is a monoid homomorphism.
It is surjective by definition and injective by uniqueness of the factorization of $M$.
Since $be=b$ for every $b\in A\cap I$, the multiplication map $(A\cap I)\times eU\to I$ is also bijective.
Both $A\cap I$ and $eU$ are submonoids of $I$ with identity $e$, and hence $I=(A\cap I)\bowtie eU$.

(2)\ The assertion follows by applying (1) to the opposite monoid $M^{\op}=U^{\op}\bowtie A^{\op}$ and its unital ideal $I^{\op}$.

Finally, if both intersections are nonempty, then $e\in A\cap U=\{1_M\}$, and hence $I=M$.
\end{proof}

\begin{cor}\label{cor:group-ideal-factors}
Let $M=A\bowtie U$ be a Zappa--Sz\'ep product of monoids, and suppose that $M$ has a group ideal $I$ with identity $e$.
Then the following hold:
\begin{enumerate}
\item If $A\cap I\neq\varnothing$, then $A\cap I$ and $U$ are both groups, and there is an exact factorization $I=(A\cap I)\bowtie eU$, where $eU\cong U$ as groups.
\item If $U\cap I\neq\varnothing$, then $A$ and $U\cap I$ are both groups, and there is an exact factorization $I=eA\bowtie(U\cap I)$, where $eA\cong A$ as groups.
\end{enumerate}
Moreover, $A\cap I\neq\varnothing$ if and only if $U$ is a group, and $U\cap I\neq\varnothing$ if and only if $A$ is a group.
In particular, $A\cap I=U\cap I=\varnothing$ if and only if neither $A$ nor $U$ is a group.
\end{cor}

\begin{proof}
(1) By Proposition~\ref{prop:unital-ideal-factors}, we have $I=(A\cap I)\bowtie eU$ and $eU\cong U$ as monoids.
Since $I$ is a group, both factors $A\cap I$ and $eU$ are groups by Remark~\ref{remark:group ZS}.
Hence $U$ is also a group.

(2) The assertion follows by applying (1) to the opposite monoid.

If $U$ is a group, choose $au\in I$ with $a\in A$ and $u\in U$.
Then $a=(au)u^{-1}\in A\cap I$.
Together with (1), this proves that $A\cap I\neq\varnothing$ if and only if $U$ is a group.
The analogous equivalence follows by symmetry.
\end{proof}

\begin{cor}\label{cor:group ideal MR}
Let $M=A\bowtie U$ be a Zappa--Sz\'ep product of monoids with a unital ideal $I$ with identity $e$.
If $A\cap I\neq\varnothing$, then there are canonical isomorphisms
\begin{align*}
\cG(M)\cong\cG(I),\quad
\cG(A)\cong\cG(A\cap I),\quad
\cG(U)\cong\cG(eU).
\end{align*}
In particular, if $I$ is a group ideal, then
\begin{align*}
\cG(M)\cong I,\quad
\cG(A)\cong A\cap I,\quad
\cG(U)\cong U.
\end{align*}
The analogous statements hold when $U\cap I\neq\varnothing$, with the roles of $A$ and $U$ interchanged.
\end{cor}

\begin{proof}
By Proposition~\ref{prop:unital-ideal-factors}, $A\cap I$ is a unital ideal of $A$ with identity $e$, and $eU\cong U$ as monoids.
The first three isomorphisms follow from Lemma~\ref{lemma:unital ideal} and the isomorphism $U\cong eU$.
If $I$ is a group ideal, then both $A\cap I$ and $U$ are groups by Corollary~\ref{cor:group-ideal-factors}, which gives the final assertion.
The case $U\cap I\neq\varnothing$ follows by passing to the opposite monoids.
\end{proof}

Suppose now that $M=A\bowtie U$ is a Zappa--Sz\'ep product of multiplicative racks.
Under the canonical isomorphisms in Corollary~\ref{cor:group ideal MR}, the exact factorization in Theorem~\ref{thm:universal-exactness} corresponds to
\begin{align*}
\cG(I)\cong\cG(A\cap I)\bowtie\cG(eU),
\end{align*}
provided that $A\cap I\neq\varnothing$.
Indeed, the retraction $\mu\colon M\to I$ restricts to the retraction $A\to A\cap I$ and the isomorphism $U\to eU$.
In particular, when $I$ is a minimal ideal, this identifies the universal group factorization with the exact factorization $I=(A\cap I)\bowtie eU$ of the minimal group ideal.

We next characterize proper unital ideals meeting a prescribed factor in terms of idempotents.

\begin{prop}\label{prop:idempotents-unital-ideals}
Let $M=A\bowtie U$ be a Zappa--Sz\'ep product of multiplicative racks.
The following are equivalent$:$
\begin{enumerate}
\item There exists a proper unital ideal $I$ of $M$ such that $A\cap I\neq\varnothing$.
\item The monoid $A$ contains an idempotent distinct from $1_M$.
\end{enumerate}
The analogous equivalence holds with $A$ replaced by $U$.
In particular, if $A$ (resp. $U$) is finite and not a group, there exists a proper unital ideal meeting $A$ (resp. $U$).
\end{prop}

\begin{proof}
Suppose that $I$ is a proper unital ideal meeting $A$, and let $e$ be its identity.
Proposition~\ref{prop:unital-ideal-factors} gives $e\in A$.
Since $I$ is proper, $e\neq1_M$, so (2) holds.
Conversely, let $e\in A$ satisfy $e^2=e\neq1_M$.
By Lemma~\ref{lemma:MRidempotent}, $e$ is central in $M$.
Hence $eM$ is a two-sided ideal with identity $e$, and it is proper because $1_M\in eM$ would imply $e=1_M$.
Since $e\in A\cap eM$, condition (1) holds.
The argument for $U$ is identical.

Finally, every element of a finite monoid has an idempotent positive power.
Thus, if the identity is its only idempotent, every element has a positive power equal to the identity and is therefore invertible.
Hence every finite non-group monoid has a nonidentity idempotent.
\end{proof}

\begin{remark}\label{rem:ideal-vs-minimal}
For a finite Zappa--Sz\'ep product $M=A\bowtie U$ of multiplicative racks with two non-group factors, neither factor meets the minimal ideal by Corollary~\ref{cor:group-ideal-factors}.
Nevertheless, Proposition~\ref{prop:idempotents-unital-ideals} shows that each factor meets some proper unital ideal, although these ideals must be distinct.
To apply Theorem~\ref{thm:MR-Bianchi}, one must also verify that the chosen ideal is invariant under all right translations.
\end{remark}

\subsection{Zappa--Sz\'ep factorizations of generalized Bianchi-type multiplicative racks}\label{subsec:BianchiZS}
Throughout this subsection, let $(Q,I,\varepsilon,\rho)$ be a PMR--MR pair and put
\begin{align*}
M\coloneq\cB(Q,I)=Q\sqcup\overline I.
\end{align*}
Let $e=1_I$, and let $\mu\colon M\to I$ be the canonical homomorphism obtained by composing the retraction $x\mapsto x\overline e$ onto $\overline I$ with the identification $\overline I\cong I$.
Thus $\mu|_Q=\varepsilon$ and $\mu(\overline i)=i$ for all $i\in I$.
We study exact factorizations of $M$ whose first factor meets the distinguished ideal $\overline I$.

We first construct multiplicative subracks of $M$ from multiplicative subracks of $I$.

\begin{lemma}\label{lem:Bianchi-restriction}
Let $B$ be a multiplicative subrack of $I$ with identity $e$, and put $P\coloneq\varepsilon^{-1}(B)$.
Then $(P,B,\varepsilon|_P,\rho|_B)$ is a PMR--MR pair.
Moreover, the natural inclusion identifies $\cB(P,B)$ with the multiplicative subrack $P\sqcup\overline B$ of $M$.
\end{lemma}

\begin{proof}
Since $\varepsilon(1_Q)=e\in B$, we have $1_Q\in P$.
If $x,y\in P$ and $xy$ is defined in $Q$, then $\varepsilon(xy)=\varepsilon(x)\varepsilon(y)\in B$, so $xy\in P$.
Since $\varepsilon$ is a rack homomorphism and $B$ is a subrack of $I$, both $S_y$ and $S_y^{-1}$ preserve $P$ for every $y\in P$.
Thus $P$ inherits a rack structure from $Q$.
The partial multiplication satisfies conditions~\eqref{PM1}--\eqref{PM5}, since products defined between elements of $P$ remain in $P$ and the domain conditions are inherited from $Q$.
Hence $P$ is a partially multiplicative rack.

For $x\in P$ and $b\in B$, condition~\eqref{GP2} gives $\varepsilon(\rho_b(x))=\varepsilon(x)*b\in B$.
Since $\rho_b$ and $S_b$ are bijective and $\varepsilon\circ\rho_b=S_b\circ\varepsilon$, we also have
\begin{align*}
\varepsilon(\rho_b^{-1}(x))=S_b^{-1}(\varepsilon(x))\in B.
\end{align*}
Thus $\rho_b$ restricts to an automorphism of $P$ as a partially multiplicative rack.
Conditions~\eqref{GP1} and \eqref{GP2} are inherited from the original PMR--MR pair, so $(P,B,\varepsilon|_P,\rho|_B)$ is a PMR--MR pair.

Finally, the multiplication and rack operation on $\cB(P,B)$ agree with the restrictions of those on $M$.
Hence the natural inclusion identifies $\cB(P,B)$ with $P\sqcup\overline B$ as a multiplicative subrack of $M$.
\end{proof}

We next show that every exact factorization whose first factor meets $\overline I$ has the form described in Lemma~\ref{lem:Bianchi-restriction}.

\begin{lemma}\label{lem:Bianchi-factor-restriction}
Suppose that $M=A\bowtie U$ is an exact factorization of multiplicative racks with $A\cap\overline I\neq\varnothing$.
Put $B \coloneq\{i\in I\mid\overline i\in A\}$ and $P \coloneq A\cap Q$.
Then $P=\varepsilon^{-1}(B)$, and the multiplication map
\begin{align*}
P\times U\to Q;\ (p,u)\mapsto pu
\end{align*}
is bijective.
\end{lemma}

\begin{proof}
By Proposition~\ref{prop:unital-ideal-factors}, we have $\overline e\in A$, an exact factorization $\overline I=(A\cap\overline I)\bowtie\overline eU$ and an isomorphism $U\to\overline eU$ given by $u\mapsto\overline e u$.
Under the identification $\overline I\cong I$, these yield an exact factorization $I=B\bowtie\mu(U)$ and the injectivity of $\mu|_U$.
Since $\overline I$ is a proper ideal and $A\cap\overline I\neq\varnothing$, the final assertion of Proposition~\ref{prop:unital-ideal-factors} gives $U\cap\overline I=\varnothing$, so $U\subseteq Q$.

Since $b\overline e=b$ for every $b\in A\cap\overline I$, we have $\overline I=(A\cap\overline I)U$.
By uniqueness of the factorization $M=A\bowtie U$, it follows that
\begin{align*}
Q=M\setminus\overline I=(A\setminus\overline I)U=PU.
\end{align*}
The same uniqueness implies that the multiplication map $P\times U\to Q$ is bijective.

It remains to show that $P=\varepsilon^{-1}(B)$.
Since $\overline e\in A$, we have $\mu(A)\subseteq B$.
For any $m=au\in M$ with $a\in A$ and $u\in U$, uniqueness of the factorization $I=B\bowtie\mu(U)$ gives
\begin{align*}
\mu(m)\in B\iff\mu(u)=e\iff u=1_M,
\end{align*}
where the last equivalence follows from the injectivity of $\mu|_U$.
Thus $\mu^{-1}(B)=A$, and hence
\begin{align*}
\varepsilon^{-1}(B)=Q\cap\mu^{-1}(B)=Q\cap A=P. &\qedhere
\end{align*}
\end{proof}

The preceding lemmas give a criterion for exact factorizations of generalized Bianchi-type multiplicative racks.

\begin{thm}\label{thm:Bianchi-ZS-factorization}
Let $B$ be a multiplicative subrack of $I$ with identity $e$, and let $U\subseteq Q$ be a multiplicative subrack of $M$.
Put $P\coloneq\varepsilon^{-1}(B)$ and identify $A\coloneq\cB(P,B)=P\sqcup\overline B$ with its canonical image in $M$.
Then $M=A\bowtie U$ is an exact factorization of multiplicative racks if and only if the following conditions hold:
\begin{enumerate}
\item The restriction $\mu|_U\colon U\to\mu(U)$ is injective.
\item The multiplicative subracks $B$ and $\mu(U)$ form an exact factorization $I=B\bowtie\mu(U)$.
\item The multiplication map $P\times U\to Q;\ (p,u)\mapsto pu$ is bijective.
\end{enumerate}
Moreover, every exact factorization $M=A\bowtie U$ of multiplicative racks with $A\cap\overline I\neq\varnothing$ arises in this way.
\end{thm}

\begin{figure}[t]
\centering
\begin{tikzpicture}[
    x=0.85cm,
    y=0.85cm,
    every node/.style={font=\small}
]
\fill[blue!15] (0,0) rectangle (6,5);
\fill[red!15] (6,0) rectangle (12,5);
\draw[thick] (0,0) rectangle (12,5);
\draw[thick] (6,0) -- (6,5);
\node[above] at (6,5.1)
{{\Large$M=\mathcal B(Q,I)=A\bowtie U$}};
\node at (3,4.1) {$Q$};
\node at (3,3.0) {$Q=PU$};
\node at (3,1.9) {%
  $\displaystyle
  P\times U\xrightarrow{\cong}Q
  \qquad\text{(3)}$
};
\node at (9,4.1) {$\overline I$};
\node at (9,3.25) {$\overline I=\overline B U$};
\node at (9,2.25) {%
  $\displaystyle
  U\xrightarrow[\text{injective}]{\mu|_U}\mu(U)
  \qquad\text{(1)}$
};
\node at (9,1.25) {%
  $\displaystyle
  I=B\bowtie\mu(U)
  \qquad\text{(2)}$
};
\draw[thick] (-0.65,0) -- (-0.65,5);
\draw[thick] (-0.85,0) -- (-0.45,0);
\draw[thick] (-0.85,5) -- (-0.45,5);
\node[left] at (-0.95,2.5) {$U$};
\draw[thick] (0,-0.85) -- (12,-0.85);
\draw[thick] (0,-1.05) -- (0,-0.65);
\draw[thick] (6,-1.05) -- (6,-0.65);
\draw[thick] (12,-1.05) -- (12,-0.65);
\node at (3,-0.4) {$P=\varepsilon^{-1}(B)$};
\node at (9,-0.4) {$\overline B$};
\node at (6,-1.35) {$A=P\sqcup\overline B$};
\end{tikzpicture}
\caption{An illustration of the factorization $M=A\bowtie U$, with $A=P\sqcup\overline B$.}
\label{fig:Bianchi-ZS}
\end{figure}
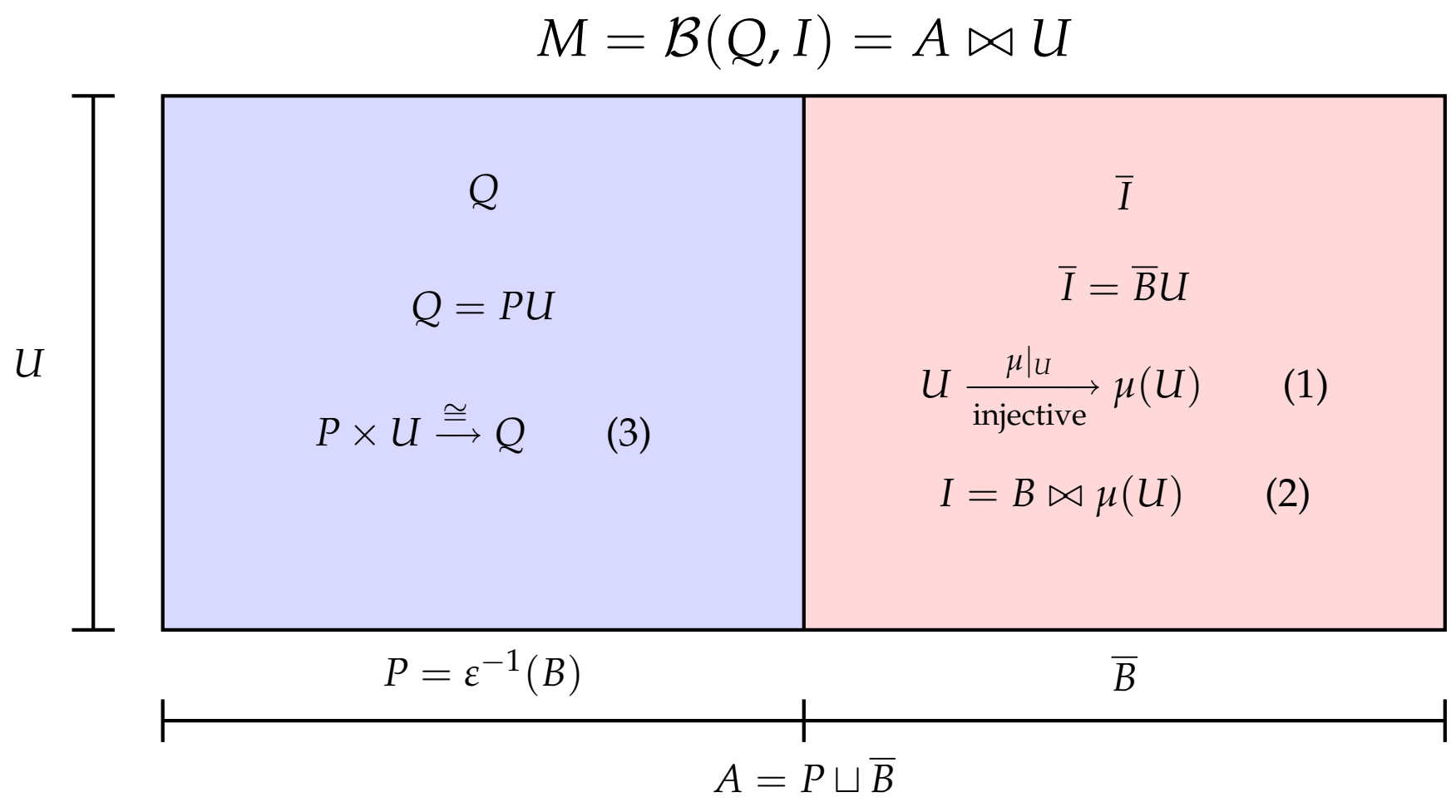

The factorization described in the theorem is illustrated in Figure~\ref{fig:Bianchi-ZS}.
The two regions correspond to the complement $Q$ and the distinguished ideal $\overline I$, with conditions (1) and (2) describing the factorization of the ideal and condition (3) describing that of the complement.

\begin{proof}
By Lemma~\ref{lem:Bianchi-restriction}, $A$ is a multiplicative subrack of $M$.
Suppose first that $M=A\bowtie U$.
Since $A\cap\overline I=\overline B\neq\varnothing$, Proposition~\ref{prop:unital-ideal-factors} gives an exact factorization $\overline I=\overline B\bowtie\overline eU$ and an isomorphism $U\to\overline eU$.
Under the identification $\overline I\cong I$, these give conditions (1) and (2), since $\mu$ is a homomorphism of multiplicative racks.
Condition (3) follows from Lemma~\ref{lem:Bianchi-factor-restriction}.

Conversely, suppose that conditions (1)--(3) hold.
By condition (3), multiplication gives a bijection $P\times U\to Q$.
By conditions (1) and (2), every $i\in I$ has a unique expression $i=b\mu(u)$ with $b\in B$ and $u\in U$.
Since $\overline b\,u=\overline{b\mu(u)}$ in $M$, multiplication also gives a bijection $\overline B\times U\to\overline I$.
Combining these bijections and the disjoint unions $A=P\sqcup\overline B$ and $M=Q\sqcup\overline I$, we obtain a bijection $A\times U\to M$.
Hence $M=A\bowtie U$.

Finally, consider an arbitrary exact factorization $M=A\bowtie U$ with $A\cap\overline I\neq\varnothing$.
By Proposition~\ref{prop:unital-ideal-factors}, we have $\overline e\in A$ and $U\subseteq Q$.
Since $A$ is a multiplicative subrack and $\overline I$ is invariant under all right translations, their intersection is a multiplicative subrack of $\overline I$.
Let $B$ be the corresponding multiplicative subrack of $I$.
By Lemma~\ref{lem:Bianchi-factor-restriction}, we have $A\cap Q=\varepsilon^{-1}(B)$.
Thus $A=\varepsilon^{-1}(B)\sqcup\overline B=\cB(\varepsilon^{-1}(B),B)$, proving the final assertion.
\end{proof}

When the second factor is a group, the factorization of the distinguished ideal automatically extends to the whole multiplicative rack.

\begin{cor}\label{cor:Bianchi-group-factor}
Let $B$ be a multiplicative subrack of $I$ with identity $e$, and let $U\leq M^\times$ be a subgroup.
Suppose that $\mu|_U$ is injective and that $I=B\bowtie\mu(U)$ is an exact factorization of multiplicative racks.
Then there is an exact factorization of multiplicative racks
\begin{align*}
M=\cB(\varepsilon^{-1}(B),B)\bowtie U,
\end{align*}
where $\cB(\varepsilon^{-1}(B),B)$ is identified with its canonical image in $M$.
\end{cor}

\begin{proof}
Since $\overline I$ is a proper ideal, it contains no units of $M$, so $U\subseteq Q$.
Put $P\coloneq\varepsilon^{-1}(B)$.
By Theorem~\ref{thm:Bianchi-ZS-factorization}, it suffices to show that the multiplication map $P\times U\to Q$ is bijective.

Take $q\in Q$.
By the factorization $I=B\bowtie\mu(U)$ and the injectivity of $\mu|_U$, there exist unique $b\in B$ and $u\in U$ such that $\varepsilon(q)=b\mu(u)$.
Put $p\coloneq qu^{-1}$.
Since $u$ is a unit and $\overline I$ is an ideal, we have $p\in Q$.
Moreover,
\begin{align*}
\varepsilon(p)=\varepsilon(q)\mu(u)^{-1}=b,
\end{align*}
so $p\in P$ and $q=pu$.
Thus $P\times U\to Q$ is surjective.

Suppose that $pu=p'u'$ for some $p,p'\in P$ and $u,u'\in U$.
Applying $\mu$ and using the uniqueness of the factorization $I=B\bowtie\mu(U)$ gives $\mu(u)=\mu(u')$.
Since $\mu|_U$ is injective, we obtain $u=u'$.
Right multiplication by a unit is injective, so $p=p'$.
Thus $P\times U\to Q$ is bijective, and the result follows from Theorem~\ref{thm:Bianchi-ZS-factorization}.
\end{proof}

\subsection{Layered constructions and non-cancellative examples}\label{subsec:Bianchi-noncancellative-examples}
We describe a layered construction of multiplicative quandles from exact factorizations of groups.
It provides a common framework for finite and infinite non-cancellative examples and allows us to compute their universal groups.

Let $H=D\bowtie C$ be an exact factorization of groups, and let $S$ and $T$ be commutative monoids.
Suppose that there are families of central subgroups
\begin{align*}
K_s\leq D\cap Z(H)\quad(s\in S),\quad J_t\leq C\cap Z(H)\quad(t\in T)
\end{align*}
such that
\begin{align*}
K_{1_S}=J_{1_T}=\{1\},\quad K_sK_{s'}\leq K_{ss'},\quad J_tJ_{t'}\leq J_{tt'}
\end{align*}
for all $s,s'\in S$ and $t,t'\in T$.

For each $(s,t)\in S\times T$, let $L_{s,t}$ be a disjoint copy of $H/(K_sJ_t)$, and denote the class of $h\in H$ in this copy by $[h]_{s,t}$.
Put
\begin{align*}
M\coloneq\bigsqcup_{(s,t)\in S\times T}L_{s,t},
\end{align*}
and define
\begin{align*}
[g]_{s,t}[h]_{s',t'}&\coloneq[gh]_{ss',tt'},\\
[g]_{s,t}*[h]_{s',t'}&\coloneq[h^{-1}gh]_{s,t}.
\end{align*}
These operations are well-defined because the subgroups $K_sJ_t$ are central and satisfy
\begin{align*}
K_sJ_t,\ K_{s'}J_{t'}\leq K_{ss'}J_{tt'}.
\end{align*}
The multiplication is associative with identity $[1_H]_{1_S,1_T}$.
Each right translation is induced by conjugation in $H$ and hence is a monoid automorphism.
The identities $S_{xy}=S_y\circ S_x$ and $xy=yS_y(x)$ follow from those in $H$ and the commutativity of $S\times T$.
Thus $M$ is a multiplicative rack by Remark~\ref{remark:MR}.
Since $x*x=x$ for every $x\in M$, it is a multiplicative quandle.

Consider the subsets
\begin{align*}
A \coloneq\bigsqcup_{s\in S}\{[d]_{s,1_T}\mid d\in D\}, \quad U \coloneq\bigsqcup_{t\in T}\{[c]_{1_S,t}\mid c\in C\}.
\end{align*}
These are multiplicative subquandles of $M$.
For each $(s,t)\in S\times T$, multiplication induces a surjective map
\begin{align*}
(D/K_s)\times(C/J_t) \to H/(K_sJ_t).
\end{align*}
This map is injective as well.
Indeed, if $d_1c_1(K_sJ_t)=d_2c_2(K_sJ_t)$, then $d_2^{-1}d_1c_1c_2^{-1}\in K_sJ_t$.
Uniqueness of the factorization $H=D\bowtie C$ gives $d_2^{-1}d_1\in K_s$ and $c_1c_2^{-1}\in J_t$.
Since the layers are disjoint, these bijections yield an exact factorization
\begin{align*}
M=A\bowtie U.
\end{align*}

The universal group of $M$ can also be computed explicitly.
Put
\begin{align*}
K\coloneq\langle K_s\mid s\in S\rangle,\quad J\coloneq\langle J_t\mid t\in T\rangle.
\end{align*}
The natural map
\begin{align*}
H\times S\times T \to M;\ (h,s,t) \mapsto[h]_{s,t}
\end{align*}
is a surjective monoid homomorphism.
Since
\begin{align*}
\cG(H\times S\times T)\cong H\times\cG(S)\times\cG(T),
\end{align*}
it induces a surjective homomorphism from this group onto $\cG(M)$.
For $k\in K_s$, the equality $[k]_{s,1_T}=[1_H]_{s,1_T}$ implies that the image of $k$ is trivial in $\cG(M)$.
The same holds for every element of $J_t$.
Thus the induced homomorphism factors through $H/(KJ)\times\cG(S)\times\cG(T)$.
Conversely, the map
\begin{align*}
M \to H/(KJ)\times\cG(S)\times\cG(T);\ [h]_{s,t} \mapsto(hKJ,\omega_S(s),\omega_T(t))
\end{align*}
is a well-defined monoid homomorphism.
The resulting group homomorphisms are mutually inverse on the generators, giving a canonical isomorphism
\begin{align*}
\cG(M)\cong H/(KJ)\times\cG(S)\times\cG(T).
\end{align*}
Similarly, we obtain
\begin{align*}
\cG(A) \cong(D/K)\times\cG(S), \quad \cG(U) \cong(C/J)\times\cG(T).
\end{align*}
Under these isomorphisms, the exact factorization of Theorem~\ref{thm:universal-exactness} corresponds to the factorization $H/(KJ)=(D/K)\bowtie(C/J)$, together with the direct-product factors $\cG(S)$ and $\cG(T)$.

We use the following elementary criterion for direct indecomposability.

\begin{lemma}\label{lem:nonunit-stabilizer-indecomposable}
Let $M$ be a monoid whose group of units is nontrivial and directly indecomposable.
Suppose that every nonunit $x\in M$ satisfies $xu=x$ for some nonidentity unit $u\in M^\times$.
Then $M$ is directly indecomposable as a monoid.
\end{lemma}

\begin{proof}
Suppose that $M\cong X\times Y$ for two nontrivial monoids $X$ and $Y$.
Since $M^\times\cong X^\times\times Y^\times$ is directly indecomposable, after interchanging the factors we may assume that $X^\times\cong M^\times$ and $Y^\times=\{1_Y\}$.
Choose $y\in Y\setminus\{1_Y\}$.
Then $(1_X,y)$ is a nonunit, but
\begin{align*}
(1_X,y)(u,1_Y)=(u,y)\neq(1_X,y)
\end{align*}
for every nonidentity $u\in X^\times$.
This contradicts the assumption.
\end{proof}

\begin{ex}\label{ex:Bianchi-unified-family}
We apply the preceding construction to obtain finite and infinite directly indecomposable multiplicative quandles with non-group factors and nontrivial Zappa--Sz\'ep actions.

Fix an odd prime $p$ and put $V\coloneq\bF_p^3$.
Define a group $H$ on $V\times\bigwedge^2V$ by
\begin{align*}
(v,w)(v',w')=(v+v',w+w'+v\wedge v').
\end{align*}
The commutator is given by
\begin{align*}
[(v,w),(v',w')]=(0,2v\wedge v').
\end{align*}
Since $p$ is odd and the elements $v\wedge v'$ span $\bigwedge^2V$, we have
\begin{align*}
[H,H]=\{0\}\times\bigwedge^2V.
\end{align*}
Moreover, $v\wedge v'=0$ for all $v'\in V$ implies $v=0$, so $Z(H)=[H,H]$.
In particular, $|H|=p^6$.

Let $e_1,e_2,e_3$ be the standard basis of $V$, and define
\begin{align*}
D&\coloneq\{(v,w)\in H\mid v\in\langle e_1,e_2\rangle,\ w\in\langle e_1\wedge e_2,e_1\wedge e_3\rangle\},\\
C&\coloneq\{(v,w)\in H\mid v\in\langle e_3\rangle,\ w\in\langle e_2\wedge e_3\rangle\}.
\end{align*}
These are subgroups of orders $p^4$ and $p^2$, respectively, with trivial intersection.
Hence $H=D\bowtie C$.
Neither subgroup is normal in $H$, since
\begin{align*}
(e_1,0)^{-1}(e_3,0)(e_1,0)&=(e_3,-2e_1\wedge e_3)\notin C,\\
(e_3,0)^{-1}(e_2,0)(e_3,0)&=(e_2,2e_2\wedge e_3)\notin D.
\end{align*}
Put
\begin{align*}
K_1\coloneq\langle(0,e_1\wedge e_3)\rangle,\quad K_2\coloneq\langle(0,e_2\wedge e_3)\rangle.
\end{align*}
These are central subgroups of order $p$ contained in $D$ and $C$, respectively, with trivial intersection.

Let $N$ be a positive integer or $\infty$.
Take
\begin{align*}
S&\coloneq\{0,1\},\\
T_N&\coloneq
\begin{cases}
\{0,1,\ldots,N\},&N<\infty,\\
\bN_0,&N=\infty,
\end{cases}
\end{align*}
with multiplication given by
\begin{align*}
s\cdot t&\coloneq\max\{s,t\}\qquad(s,t\in S),\\
n\cdot m&\coloneq
\begin{cases}
\min\{n+m,N\},&N<\infty,\\
n+m,&N=\infty.
\end{cases}
\end{align*}
Both monoids have identity $0$.
Define families of central subgroups by
\begin{align*}
K_0^S \coloneq\{1\}, \quad K_1^S \coloneq K_1, \quad J_0^T \coloneq\{1\}, \quad J_n^T \coloneq K_2\quad(n\geq1).
\end{align*}
These families satisfy the assumptions of the layered construction.
Let $M_N=A\bowtie U_N$ be the resulting multiplicative quandle and its exact factorization.

Both factors are non-groups, since $A$ contains a nonidentity idempotent and $U_N$ contains noninvertible elements of positive layer index.
Moreover, both components of the Zappa--Sz\'ep action are nontrivial, since the factorization restricts to $H=D\bowtie C$ on the identity layer and neither $D$ nor $C$ is normal in $H$.

Put
\begin{align*}
I_N\coloneq\bigsqcup_{n\in T_N}L_{1,n},\quad e\coloneq[1]_{1,0}.
\end{align*}
Then $I_N=eM_N$ is a proper unital ideal with identity $e$, but it is not a group because it contains elements of positive layer index.
Every right translation preserves each layer, so Theorem~\ref{thm:MR-Bianchi} gives a generalized Bianchi-type presentation $M_N\cong\cB(Q_N,I_N)$, where $Q_N=M_N\setminus I_N$.
Since $A\cap I_N\neq\varnothing$, Theorem~\ref{thm:Bianchi-ZS-factorization} shows that the factorization $M_N=A\bowtie U_N$ arises from this construction.

We show that every $M_N$ is directly indecomposable.
First, the group $H$ is directly indecomposable.
Suppose that $H=H_1\times H_2$ is a nontrivial direct product.
Then $V\cong H/Z(H)$ decomposes as $V=V_1\oplus V_2$, where $V_i\cong H_i/Z(H_i)$.
Since the two factors commute, we have $v_1\wedge v_2=0$ for all $v_i\in V_i$.
If both $V_1$ and $V_2$ are nonzero, choosing nonzero $v_i\in V_i$ gives linearly independent vectors with $v_1\wedge v_2\neq0$, a contradiction.
Thus one factor is abelian, which must be trivial because $Z(H)=[H,H]$.
Hence $H$ is directly indecomposable.

Since $(0,0)$ is the only unit of the indexing monoid $S\times T_N$, the group of units of $M_N$ is $M_N^\times=L_{0,0}\cong H$.
Every nonunit belongs either to $L_{1,0}$ or to a layer $L_{s,n}$ with $n\geq1$.
Elements of $L_{1,0}$ are fixed by right multiplication by $K_1$, whereas elements of $L_{s,n}$ with $n\geq1$ are fixed by right multiplication by $K_2$.
Since both subgroups are nontrivial and $H$ is directly indecomposable, Lemma~\ref{lem:nonunit-stabilizer-indecomposable} shows that $M_N$ is directly indecomposable as a monoid.

Suppose first that $N<\infty$.
Then
\begin{align*}
|M_N| =p^6+p^5+N(p^5+p^4), \quad |A| =p^4+p^3, \quad |U_N| =p^2+Np.
\end{align*}
Since $(1,N)$ is the zero of the indexing monoid $S\times T_N$, the layer $L_{1,N}\cong H/(K_1K_2)$ is a group ideal of $M_N$.
Hence it is the unique minimal ideal of $M_N$ of order $p^4$.
For $N\geq2$, the element $[1]_{0,1}$ is not regular.
Indeed, the second layer index of $[1]_{0,1}y[1]_{0,1}$ is at least $2$ for every $y\in M_N$, whereas that of $[1]_{0,1}$ is $1$.

On the other hand, $M_\infty$ is countably infinite and has no minimal ideal.
Let $J$ be any ideal of $M_\infty$, and choose $x\in J\cap L_{s,n}$.
The subset
\begin{align*}
T_{n+1}'\coloneq\bigsqcup_{\substack{t\in\{0,1\}\\m\geq n+1}}L_{t,m}
\end{align*}
is an ideal of $M_\infty$.
Since $x[1]_{0,1}\in J\cap T_{n+1}'$ but $x\notin T_{n+1}'$, the intersection $J\cap T_{n+1}'$ is a proper ideal contained in $J$.
Thus $M_\infty$ has no minimal ideal.

Finally, the universal group formula gives
\begin{align*}
\cG(M_N)\cong
\begin{cases}
H/(K_1K_2),&N<\infty,\\
H/(K_1K_2)\times\bZ,&N=\infty.
\end{cases}
\end{align*}
Thus the same construction produces finite and countably infinite directly indecomposable non-cancellative multiplicative quandles with non-group factors and nontrivial Zappa--Sz\'ep actions.
\end{ex}


\bibliographystyle{amsalpha}
\bibliography{MR}
\end{document}